\documentclass{amsart}
\usepackage{amssymb}
\usepackage{stmaryrd}
\usepackage[all]{xy}

\title{CYLINDRIC ALGEBRAS OF DE MORGAN-VALUED LOGIC}

\author{}

\address{Mathematics Department\\Sonoma State University\\Rohnert Park, CA 94928, USA}

\email{norm.feldman@sonoma.edu}

\newtheorem{Theorem}{Theorem}[section] 

\newtheorem{Lemma}[Theorem]{Lemma}

\theoremstyle{definition}
\newtheorem{Definition}[Theorem]{Definition}

\begin{document}

\begin{abstract}
We construct a De Morgan algebra-valued logic with quantifiers, where the truth values are in a finite De Morgan algebra, We show that there is a representation theorem of the cylindric algebra of this logic from which a completeness theorem for De Morgan algebra-valued logic follows. This is a generalization of the results in [2].
\end{abstract}
\maketitle

\section{Introduction}
We present a set of axioms for generalized cylindric algebras based on a finite De Morgan algebra $\mathfrak{M}$. These cylindric algebras are referred to as $\mathfrak{M}$-cylindric algebras and will be used to define the semantics of a logic with quantifiers with truth values in $\mathfrak{M}$. We show how to construct an $\mathfrak{M}$-cylindric algebra from a cylindric algebra. A representation theorem for $\mathfrak{M}$-cylindric algebras is proved. The completeness theorem for this logic follows from the representation theorem. This result generalizes the completeness theorem proved in [2].

\section{Cylindric Algebra of De Morgan-Valued Logic}

\subsection{De Morgan Algebras}

\begin{Definition} A De Morgan algebra is an algebra $\langle M, +, \cdot, -, 0, 1\rangle$ where the following are satisfied:

\begin{tabbing}
xxxxxxxxxxxxxxxxxxxxxxxxxxxxxxx\= \kill

1 (a) $x+y=y+x$ \> (b) $x\cdot y = y\cdot x$\\

2 (a) $(x+y)+z=x+(y+z)$ \> (b) $(x\cdot y)\cdot z=x\cdot (y\cdot z)$\\

3 (a) $x \cdot (y+z)=(x\cdot y)+(x\cdot z)$ \> (b) $x + (y\cdot z)=(x+ y)\cdot(x+ z)$\\

4 (a) $x+0=x$ \> (b) $x\cdot 1=x$\\

5 (a) $x\cdot (x+y)=x$ \> (b) $x+ (x\cdot y)=x$\\

6 (a) $-(x+y)=-x\cdot -y$ \> (b) $-(x\cdot y)=-x+ -y$\\

7 $--x=x$\\

\end{tabbing}

\end{Definition}

The following are easily verified:

\begin{tabbing}
xxxxxxxxxxxxxxxxxxxxxxxxxxxxxxx\= \kill
1 (a) $-0=1$  \> (b) $-1=0$\\

2 (a) $x+ x=x$  \> (b) $x\cdot x=x$\\

3 (a) $0 \cdot x=0$  \> (b) $x+1=1$\\

\end{tabbing}

In the following we will assume that  $\mathfrak{M}=\langle M, +, \cdot, -, 0, 1\rangle$ is a finite De Morgan algebra with $n$ elements. A binary relation on $M$ is defined as follows: $x\leq y$ if $x\cdot y=x$, or equivalently, $x+y=y$. It can be shown that $\leq$ is a partial ordering on $M$.

We refer to the cylindric algebra associated with the De Morgan algebra $\mathfrak{M}$ as an $\mathfrak{M}$-cylindric algebra.

\subsection{ Axioms for $\mathfrak{M}$-Cylindric Algebras}
 
\begin{Definition} 
In the following let $\alpha$ be an ordinal. 
$\mathfrak{A}=\langle A,+,\cdot,-,u_p,c_\kappa,d_{\kappa \lambda},\delta_p \rangle_{\kappa,\lambda <\alpha,p \in M}$ is a $\mathfrak{M}$-cylindric algebra of dimension $\alpha$ ($\mathfrak{M}-\textrm{CA}_{\alpha}$) if $+$ and $\cdot$ are binary operations on $A$; $-$, $C_\kappa$, and $\delta_p$ for $p \in M$ and $\kappa<\alpha$ are unary operations on $A$; and $u_p$ and $d_{\kappa \lambda}$ are members of $A$ for $p\in M$ and $\kappa,\lambda< \alpha$ satisfying the following axioms:

Boolean Axioms:

\begin{tabbing}
xxxxxxxxxxxxxxxxxxxxxxxxxxxxxxxxxxxxxx\= \kill
1 (a) $a+b=b+a$ \> (b) $a \cdot b=b \cdot a$\\
2 (a) $(a+b)+c=a+(b+c)$ \> (b) $(a \cdot b) \cdot c=a \cdot (b \cdot c)$\\
3 (a) $a \cdot (b+c)=(a \cdot b)+(a \cdot c)$ \>(b)  $a+(b \cdot c)=(a+b) \cdot (a+c)$\\
4 (a) $a+u_0=a$ \> (b) $a \cdot u_1=a$\\
5 (a) $\delta_pa+-\delta_pa=u_1$ $(p\in M)$\> (b) $\delta_p a \cdot -\delta_p a=u_0$ $(p\in M)$\\
\end{tabbing}

z
Cylindric Axioms:

\begin{tabbing}
xxxxxxxxxxxxxxxxxxxxxxxxxxxxxxxxxxxxxx\= \kill
6 $c_\kappa u_0=u_0$  \> 7 $a+c_\kappa a=c_\kappa a$\\
8 $c_\kappa(a \cdot c_\kappa b)=c_\kappa a \cdot c_\kappa b$ \> 9 $c_\kappa c_\lambda a=c_\lambda c_\kappa a$\\
10 $d_{\lambda \mu}=c_\kappa(d_{\lambda\kappa} \cdot d_{\kappa \mu})$ \> 11 $d_{\kappa \kappa}=u_1$\\
12 $c_\kappa(d_{\kappa \lambda} \cdot \delta_1a) \cdot c_\kappa(d_{\kappa \lambda} \cdot -\delta_1 a)=u_0$ $(\kappa \ne \lambda)$\\
\end{tabbing}

Since $+$ and $\cdot$ are both commutative and associative, we may use $\Sigma$ and $\prod$ for general summation and product, respectively.

\begin{tabbing}
xxxxxxxxxxxxxxxxxxxxxxxxxxxxxxxxxxxxxx\= \kill
13 $\delta_1(\delta_1a+\delta_1b)=\delta_1a+\delta_1b$  \>  14 $\delta_1(\delta_1a\cdot \delta_1b)=\delta_1a\cdot \delta_1b$\\
15 $\delta_1-\delta_1a=-\delta_1a$  \>  16 $\delta_1c_\kappa\delta_1a=c_\kappa\delta_1a$\\
17 $\delta_1d_{\kappa \lambda}=d_{\kappa \lambda}$\>  18 $\delta_0d_{\kappa \lambda}=-d_{\kappa \lambda}$\\
19 $\delta_pd_{\kappa \lambda}=u_0$ $(p \in M$ and $p\ne 0, 1)$  \>  20 $\delta_p(-a)=\delta_{-p}a$ $(p\in M)$\\
21 $ \delta_pu_p=u_1$ $(p \in M)$   \>  22 $ \delta_pu_q=u_0$ $(p,q \in M$ and $p \ne q)$\\
23 $ \delta_pa \cdot  \delta_qa=u_0$ $(p,q \in M, p\ne q)$  \>  24 $ \sum_{p \in M} \delta_p a=u_1$\\
25 $\delta_p(a+b)=\sum_{q+r=p}\delta_qa \cdot \delta_rb$  $(p \in M)$  \>  26 $\delta_p(a\cdot b)=\sum_{q\cdot r=p}\delta_qa \cdot \delta_rb$  $(p \in M)$\\
27 $\delta_1\delta_qa=\delta_qa$  \>  28 $\delta_p\delta_qa=u_0$ $(p \ne 0,1)$\\
29 $\delta_0\delta_qa=-\delta_qa$  \>  30 If $\delta_pa=\delta_pb$ for all $p \in M$, then $a=b$\\
31 $\delta_pc_\kappa a=\sum_{A \subseteq M, sup A=p}\prod_{q \in M}c_\kappa\delta_qa \cdot -\sum_{A \subseteq M, sup A>p}\prod_{q \in M}c_\kappa\delta_qa$ $(p\in M)$\\

\end{tabbing}
\end{Definition}
\begin{Theorem}
In an $\mathfrak{M}-\textrm{CA}_{\alpha}$, $a+a=a$
\begin{proof}
First we show that $\delta_pa\cdot\delta_pa=\delta_pa$. Using Axioms 3, 4, and 5 we obtain $\delta_pa\cdot\delta_pa=(\delta_pa\cdot\delta_pa)+u_0=(\delta_pa\cdot\delta_pa)+(\delta_pa\cdot-\delta_pa)=\delta_pa\cdot(\delta_pa+-\delta_pa)=\delta_pa\cdot u_1=\delta_pa$. Using Axioms 25 and 23 we obtain $\delta_p(a+a)=\sum_{q+r=p}\delta_qa\cdot\delta_ra=\delta_pa\cdot\delta_pa=\delta_pa$. Since $\delta_p(a+a)=\delta_pa$ for all $p\in M$, by Axiom 30, $a+a=a$.
\end{proof}
\end{Theorem}
\begin{Definition}
In an $\mathfrak{M}-\textrm{CA}_{\alpha}$, define $a\le b$ if $a+b=b$.
\end{Definition}

\begin{Theorem}
$\le$ is a partial order.
\end{Theorem}

\subsection{ Construction of an $\mathfrak{M}$-Cylindric Algebras}
In the following, $B^M$ is the set of functions with domain $M$ and range a subset of $B$. We write $x^p$ for $x(p)$. If $\mathfrak{A}$ is an algebra, then we denote its universe by $\vert\mathfrak{A}\vert$. Also, $\alpha$ is an ordinal and $\mathfrak{B}=\langle B,+,\cdot,-,0,1,c_\kappa,d_{\kappa \lambda} \rangle_{\kappa , \lambda < \alpha}$ is a cylindric algebra of dimension $\alpha$ ($\textrm{CA}_{\alpha}$) [2]. $c_\kappa^\partial x=-c_\kappa-x$. Note that $c_\kappa^\partial x\le x\le c_\kappa x$.

We construct an $\mathfrak{M}-\textrm{CA}_{\alpha}$ $\mathfrak{M(B)}=\langle F,\boldsymbol{+},\boldsymbol{\cdot},\boldsymbol{-},\boldsymbol{u}_p,\boldsymbol{c}_\kappa,\boldsymbol{d}_{\kappa \lambda},\boldsymbol{\delta}_p\rangle_{\kappa , \lambda < \alpha,p \in M}$ as follows:

Let $F=\{x:x \in B^M$ and $x^p\cdot x^q=0$ if $p\ne q$ and $\sum_{p\in M}x^p=1\}$.
Define $\boldsymbol{+}$ on $F$ by $(x\boldsymbol{+}y)^p=\sum_{q+r=p}x^q\cdot y^r$ for $p\in M$ and $x$, $y \in F$. To show that $x\boldsymbol{+}y \in F$, assume that $p \ne s$.  The following computation takes place in $\mathfrak{B}$. $(x\boldsymbol{+}y)^p \cdot (x\boldsymbol{+}y)^s=\sum_{q+r=p}x^q\cdot y^r\cdot \sum_{q_1+r_1=s}x^{q_1}\cdot y^{r_1}=\sum_{q_1+r_1=s}(\sum_{q+r=p}x^q\cdot y^r) \cdot x^{q_1}\cdot y^{r_1}=\sum_{q_1+r_1=s}\sum_{q+r=p}x^q\cdot y^r\cdot x^{q_1}\cdot y^{r_1}=0$ since any sumand with $q\ne q_1$ or $r_1\ne r_1$ is $0$ and if there are any remaining summands, we would have $q_1+r_1=s$, $q+r=p$, $q_1=q$, and $r_1=r$. and hence, $s=p$ contradicting $s\ne p$. In addition, $\sum_{p\in M}(x+y)^p=\sum_{p\in M}\sum_{q+r=p}x^q\cdot y^r=\sum_{q,r\in M}x^q\cdot y^r=\sum_{r\in M} \sum_{q\in M}x^q\cdot y^r=\sum_{r\in M}(y^r\cdot \sum_{r\in M}x^q)=\sum_{r\in M}y^r\cdot 1=\sum_{r\in M}y^r=1$. Therefore, $F$ is closed under $\boldsymbol{+}$.

Define $\boldsymbol{\cdot}$ on $F$ by $(x\boldsymbol{\cdot}y)^p=\sum_{q\cdot r=p}x^q\cdot y^r$ for $x$, $y \in F$. Similar to above, $F$ is closed under $\boldsymbol{\cdot}$.

Define $\boldsymbol{-}$ on $F$ by $(\boldsymbol{-} x)^p=x^{-p}$. It is easily shown that  $F$ is closed under $\boldsymbol{-}$.

Define $\boldsymbol{u}_q$ for $q \in M$ as follows: $(\boldsymbol{u}_q)^p=1$ if $q=p$ and $(\boldsymbol{u}_q)^p=0$ if $q\ne p$. Clearly, $\boldsymbol{u}_q \in F$ for all $q\in M$.

Note that throughout, whenever we write $supA$ we asssume that $A\ne\emptyset$ and $A\subseteq M$.

Define $\boldsymbol{c}_\kappa$ on $F$ by 
\begin{center}
$(\boldsymbol{c}_\kappa x)^p=\sum_{A\subseteq M,supA=p}\prod_{q\in A}c_\kappa x^q\cdot -\sum_{A\subseteq M,supA>p}\prod_{q\in A}c_\kappa x^q$
\end{center}
for $in M$.

This definition will be motivated in section 3.2.

Define $Z_\kappa(A,x)=\prod_{p\in A}c_\kappa x^p$, where $x\in B$. Let $p\ne q$. We show that $(\boldsymbol{c}_\kappa x)^p\cdot (\boldsymbol{c}_\kappa x)^q=0$ if $p\ne q$. Note that $Z_\kappa(A,x)\cdot Z_\kappa(B,x)=Z_\kappa(A\cup B,x)$

$(\boldsymbol{c}_\kappa x)^p\cdot (\boldsymbol{c}_\kappa x)^q$

$=\sum_{supA=p}Z_\kappa(A,x)\cdot -\sum_{supA>p}Z_\kappa(A,x)\cdot \sum_{supA=q}Z_\kappa(A,x)\cdot -\sum_{supA>q}Z_\kappa(A,x)$

$=\sum_{supA=p}Z_\kappa(A,x)\cdot \sum_{supA=q}Z_\kappa(A,x) \cdot \prod_{supA>p}-Z_\kappa(A,x)\cdot \prod_{supA>q}-Z_{\kappa}(A,x)$

$=\sum_{supA=p}\sum_{supB=q}[Z_\kappa(A,x)\cdot Z_\kappa(B,x)\cdot \prod_{supC>p}-Z_\kappa(C,x)\cdot \prod_{supD>q}-Z_{\kappa}(D,x)]$

Let $sup A=p$ and $sup B=q$. Clearly, $sup(A \cup B) \ge sup A=p$ and $sup(A \cup B) \ge sup B=q$. Since $p\ne q$, $sup(A\cup B)\ne p$ or $sup(A\cup B)\ne q$. Therefore, $sup(A \cup B)>p$ or $sup(A \cup B)>q$. Say, $sup(A \cup B)>p$. We then have\\ $Z_\kappa(A,x)\cdot Z_\kappa(B,x)\cdot \prod_{supC>p}-Z_\kappa(C,x)\cdot  \prod_{supD>q}-Z_\kappa(D,x)\\\le Z_\kappa(A,x)\cdot Z_\kappa(B,,x)\cdot -Z_\kappa(A\cup B,x)\\
= Z_\kappa(A\cup B,x)\cdot -Z_\kappa(A\cup B,x)=0$\\
Therefore, $(\boldsymbol{c}_\kappa x)^p\cdot (\boldsymbol{c}_\kappa x)^q=0$. Define $V_n\subseteq M$ for $n \in \omega$ as follows:

We now show that $\sum_{p\in M}(\boldsymbol{c}_\kappa x)^p=1$. 

Define $q\lessdot p$ for $q,p\in M$  by $q\lessdot p$ if $q<p$ and there is no $r\in M$ such that $p<r<q$. 

Define $V_n\subseteq M$ for $n \in \omega$ as follows: Let $V_1=\{1\}$ and $V_{n+1}=\{p\in M:p\lessdot r$ for some $r\in V_n \}$.

 Since $M$ is finite, for all $a,b\in M$ if $a<b$, there is a chain $a\lessdot c_{1}\lessdot c_{2}\lessdot \cdots \lessdot c_n \lessdot b$ and if $a\in V_k$, then $c_i\in V_{k-i}$ and $b\in V_{k-(n+1)}$ and hence, for some $m<k$ , $b\in V_m$. So we have:
 \begin{center}
 If $a<b$ and $a\in V_k$, then for some $m<k$, $b\in V_m$.   (*)
 \end{center}

By induction on n, we show that  for all $s\in V_n$, 

\begin{center}
$\sum_{supA=s}Z_\kappa(A,x)\le \sum_{p\in M}(\boldsymbol{c}_\kappa x)^p$ (**)

\end{center}

n=1: Let $s\in V_1$. Then $s=1$ and we have

\begin{center}
$\sum_{supA=s}Z_\kappa(A,x)=\sum_{supA=1}Z_\kappa(A,x)=(\boldsymbol{c}_\kappa x)^1\le \sum_{p\in M}(\boldsymbol{c}_\kappa x)^p$

\end{center}

Assume that for all $s\in V_k, \sum_{supA=s}Z_\kappa(A,x)\le  \sum_{p\in M}(\boldsymbol{c}_\kappa x)^p$ for all $k<n$. Let $s\in V_n$ and let $B=\{r\in M:r>s\}$. Then 

\begin{center}
$\{Z_\kappa(A,x):supA>s\}=\{Z_\kappa(A,x):supA=r$ for some $r\in B\}$.
\end{center}

Therefore, $\sum_{supA>s}Z_\kappa(A,x)=\sum_{r\in B}\sum_{supA=r}Z_\kappa (A,x)$. Let $r \in B$. Then $r>s$ and hence by (*), $r\in V_k$ for some $k<n$. By the induction hypothesis, $\sum_{supA=r}Z_\kappa (A,x)\le \sum_{p\in M}(\boldsymbol{c}_\kappa x)^p$. Therefore, $\sum_{supA>s}Z_\kappa(A,x)=\sum_{r\in B}\sum_{supA=r}Z_\kappa (A,x)\le \sum_{p\in M}(\boldsymbol{c}_\kappa x)^p$.

Since $(\boldsymbol{c}_\kappa x)^s\le \sum_{p\in M}(\boldsymbol{c}_\kappa x)^p$, we have $(\boldsymbol{c}_\kappa x)^s+\sum_{supA>s}Z_\kappa(A,x)\le \sum_{p\in M}(\boldsymbol{c}_\kappa x)^p$. Therefore, $\sum_{supA=s}Z_\kappa(A,x)\le \sum_{supA=s}Z_\kappa(A,x)+\sum_{supA>s}Z_\kappa(A,x)$

$=(\sum_{supA=s}Z_\kappa(A,x)\cdot -\sum_{supA>s}Z_\kappa(A,x))+\sum_{supA>s}Z_\kappa(A,x)$

$=(\boldsymbol{c}_\kappa x)^s+\sum_{supA>s}Z_\kappa(A,x)$

$\le \sum_{p\in M}(\boldsymbol{c}_\kappa x)^p$

Therefore, $\sum_{supA=s}Z_\kappa(A,x)\le \sum_{p\in M}(\boldsymbol{c}_\kappa x)^p$ for all $s\in V_n$ which proves (**).

Since $sup\{s\}=s$, $x^s\le c_\kappa x^s=Z_\kappa(\{s\},x)\le \sum_{supA=s}Z_\kappa(A,x)$. So, for all $s \in M$, $x^s\le \sum_{supA=s}Z_\kappa(A,x)\le \sum_{p\in M}(\boldsymbol{c}_\kappa x)^p$. Hence, $1=\sum_{s\in M}x^s\le \sum_{p\in M}(\boldsymbol{c}_\kappa x)^p$ and we have $\sum_{p\in M}(\boldsymbol{c}_\kappa x)^p=1$.

Therefore, $F$ is closed under $\boldsymbol{c}_\kappa$.

Define $\boldsymbol{d}_{\kappa \lambda}\in F$ by 

\begin{displaymath}
(\boldsymbol{d}_{\kappa \lambda})^p=\left \{ \begin{array}{ll}
d_{\kappa \lambda} & \textrm{if $p=1$}\\
0 &  \textrm{if $p\ne 0,1$}\\
-d_{\kappa \lambda} & \textrm{if $p=0$}\end{array} \right .
\end{displaymath}

for $p\in M$.

Clearly, $\boldsymbol{d}_{\kappa \lambda} \in F$.

Define $\boldsymbol{\delta_q}$ on $F$ by

\begin{displaymath}
(\boldsymbol{\delta}_{q}x)^p=\left \{ \begin{array}{ll}
x^q & \textrm{if $p=1$}\\
0 &  \textrm{if $p\ne 0,1$}\\
-x^q & \textrm{if $p=0$}\end{array} \right .
\end{displaymath}

for $q\in M$ 

Clearly $F$ is closed under $\delta_q$.

The above shows the following:

\begin{Theorem}
\label{is algebra}
If $\mathfrak{B}$ is a $\textrm{CA}_{\alpha}$, then $\mathfrak{M(B)}$ is an algebra.
\end{Theorem}

\begin{Theorem}
\label{less than}
In $\mathfrak{M(B)}$, $a\boldsymbol{+}b=b$ iff $b^p\le \sum_{q\le p}a^q$ for all $p\in M$.
\begin{proof}

Let  $a \boldsymbol{+} b=b$. Therefore, $(a\boldsymbol{+}b)^p=b^p$, for all $p\in M$ and by definition of $\boldsymbol{+}$ in $\mathfrak{M(B)}$, $\sum_{q+r=p}a^q\cdot b^r=b^p$ for all $p\in M$. For all $q,r\in M$ with $q+r=p$ we have $q\le p$ and $a^q\cdot b^r \le a^q$. Therefore, $b^p= \sum_{q+r= p}a^q\cdot b^r\le \sum_{q\le p}a^q$.

Now let $b^p\le \sum_{q\le p}a^q$ for all $p\in M$. Assume that $r+q=p$. We show that 
\begin{center}
$a^r\cdot b^q\le b^p$ (*).
\end{center}
If $a^r\cdot b^q=0$, then we have (*). So, assume that $a^r\cdot b^q\ne 0$.Then
\begin{center}
$0\ne a^r\cdot b^q\le a^r\cdot \sum_{t\le q}a^t=\sum_{t\leq q}a^r\cdot a^t$
\end{center}
Therefore, $\sum_{t\leq q}a^r\cdot b^t\ne 0$. However, $a^r\cdot a^t=0$ if $t\ne r$. So we must have $r=t$ for some $t\le q$. Therefore, $r
\le q$ and hence, $q=r+q=p$ and we have $a^r\cdot b^q\le b^q=b^p$ and we verified (*).

By (*), $(a\boldsymbol{+}b)^p=\sum_{r+q=p}a^r\cdot b^q \le b^p$, and hence, $(a+b)^p\le b^p$.

Since, $b^p\le \sum_{q\le p}a^q$,

$b^p=b^p\cdot b^p\le b^p\cdot \sum_{q\le p}a^q=\sum_{q\le p}b^p\cdot a^q=\sum_{p+q=p}b^p\cdot a^q \le \sum_{r+q=p}b^r\cdot a^q=(a\boldsymbol{+}b)^p$.

Here, $p$ is fixed and the sums are on $q$ and $r$. Therefore, $b^p\le(a\boldsymbol{+}b)^p$ for all $p\in M$. Since $(a\boldsymbol{+}b)^p\le b^p$ we have $b^p=(a\boldsymbol{+}b)^p$ for all $p\in M$ and hence, $b=a \boldsymbol{+}b$.
\end{proof}
\end{Theorem}

\begin{Theorem}
\label{insert}
In $\mathfrak{M(B)}$, let 
\begin{displaymath}
A^*=\{a\in \vert \mathfrak{M(B)}\vert:a^p=0 \textrm{ for all } p \in M, p\ne 0,1\}
\end{displaymath}

Let $a\boldsymbol{\mapsto} b=-a\boldsymbol{+} b$, $a\boldsymbol{\gets \mapsto} b=(a\boldsymbol{\mapsto} b)\boldsymbol{\cdot}(b\boldsymbol{\mapsto} a)$, and $\boldsymbol{c}_\kappa^\partial x=\boldsymbol{-}\boldsymbol{c}_\kappa \boldsymbol{-}x$ in $\mathfrak{M(B)}$ and $c_\kappa^\partial x=-c_\kappa -x$ in $ \mathfrak{B}$. Then\\
1. $A^*$ is closed under $\boldsymbol{+}$, $\boldsymbol{\cdot}$, $\boldsymbol{-}$, $\boldsymbol{c}_\kappa$, $\boldsymbol{c}_\kappa^\partial$, $\boldsymbol{\mapsto}$, and $\boldsymbol{\gets \mapsto}$\\
For all $a,b \in A^*$: \\
2. $a^0=-a^1$ and $a^1=-a^0$\\
3. $\boldsymbol{\delta}_pa\in A^*$ for all $a\in \vert \mathfrak{M(B)}\vert$ and $p\in M$\\
4.  $\boldsymbol{d}_{\kappa \lambda} \in A^*$ for all $\kappa, \lambda < \alpha$\\
5. $a=b$ iff $a^1=b^1$ and $a=b$ iff $a^0=b^0$\\
\begin{tabbing}
xxxxxxxxxxxxxxxxxxxxxxxxxxxxxxxxxx\= \kill
6(a) $(a\boldsymbol{+}b)^1=a^1+b^1$\>(b) $(a\boldsymbol{+}b)^0=a^0\cdot b^0$\\
7(a) $(a\boldsymbol{\cdot}b)^1=a^1\cdot b^1$\>(b) $(a\boldsymbol{\cdot}b)^0=a^0+ b^0$\\
8(a) $(\boldsymbol{-}a)^1=-a^1$\>(b) $(\boldsymbol{-}a)^0=-a^0$\\
9(a) $(\boldsymbol{c}_\kappa a)^1=c_\kappa a^1$\>(b) $(\boldsymbol{c}_\kappa a)^0=\boldsymbol{c}_\kappa^\partial a^0$\\
10(a) $(\boldsymbol{c}_\kappa^\partial a)^1=\boldsymbol{c}_\kappa^\partial a^1$ \> (b) $(\boldsymbol{c}_\kappa^\partial a)^0=\boldsymbol{c}_\kappa a^0$\\
11 $(a\boldsymbol{\mapsto} b)^1=-a^1+b^1$\\
12 $(a\boldsymbol{\gets \mapsto}b)^1=-(a^1\cdot-b^1)+(b^1\cdot-a^1)$
\end{tabbing}

\begin{proof}
As an example we will prove 9 (b). The proofs of the rest follow easily from the definitions.\\
\begin{eqnarray*}
(\boldsymbol{c}_\kappa a)^0 &=& \sum_{supA=0}\prod_{q\in A}c_\kappa a^q\cdot -\sum_{supA>0}\prod_{q\in A}c_\kappa a^q\\
&=& c_\kappa a^0\cdot -\sum_{supA>0}\prod_{q\in A}c_\kappa a^q\\
&=& c_\kappa a^0\cdot -(c_\kappa a^0\cdot c_\kappa a^1+c_\kappa a^1)\\
&=& c_\kappa a^0\cdot -c_\kappa a^1\\
&=& c_\kappa a^0\cdot -c_\kappa -a^0\\
&=& c_\kappa a^0\cdot c^\partial_\kappa a^0\\
&=& c^\partial_\kappa a^0 (\textrm{ since in a CA}_\alpha,  c^\partial_\kappa x\le c_kx)
\end{eqnarray*}
The third equality follows from the fact that if $A$ has elements other than $0$ or $1$, say $r$, then $a^r=0$ and $c_\kappa a^r=0$.
\end{proof}
\end{Theorem}

We will need the following two Boolean identities to establish that if $\mathfrak{B}$ is a $\textrm{CA}_{\alpha}$, then $\mathfrak{M(B)}$ is an $\mathfrak{M}-\textrm{CA}_{\alpha}$. We thank Dr. Brannen of the Mathematics Department at Sonoma State University for his help in discovering and proving the two identities.\\
\begin{Theorem}
\label{boolean 1}
Let $\mathfrak{B}=\langle B,+,\cdot ,-,0,1\rangle$ be a Boolean algebra, $p\in M$,\\$U,W:M\to B$, $\sum_{q\in M}W_q=1$, $W_q\cdot W_r=0$ if $q\ne r$. Define $V:M\to B$ by $V_r=\sum_{supA=r}\prod_{q\in A}U_q\cdot -\sum_{supA>r}\prod_{q\in A}U_q$ for $r \in M$. Then for all $p\in M$\\
\begin{center}
$\sum_{supA=p}\prod_{q\in A}\sum_{r\cdot t=q}(U_r\cdot W_t)\cdot-\sum_{supA>p}\prod_{q\in A}\sum_{r\cdot t=q}(U_r\cdot W_t)=\sum_{r\cdot t=p}(V_r\cdot W_t)$
\end{center}
\begin{proof}
Since every Boolean algebra is isomorphic to a Boolean algebra of set, we will prove the following:
Let $Z$ be a subset of the set of subsets of a set $X$ such that $Z$ is closed under $\cup$, $\cap$, $-$ (complementation with respect to $X$), and $\emptyset$, $X \in Z$, $p\in M$, $U,W:M\to Z$, $\bigcup_{q\in M}W_q=X$, $W_q\cap W_r=\emptyset$ if $q\ne r$ and $V:M\to Z$ is defined by $V_r=\bigcup_{supA=r}\bigcap_{q\in A}U_q\cap -\bigcup_{supA>r}\bigcap_{q\in A}U_q$ for $r \in M$. Then for all $p\in M$\\
\begin{center}
$\bigcup_{supA=p}\bigcap_{q\in A}\bigcup_{r\cdot t=q}(U_r\cap W_t)\cap -\bigcup_{supA>p}\bigcap_{q\in A}\bigcup_{r\cdot t=q}(U_r\cap W_t) =\bigcup_{r\cdot t=p}(V_r\cap W_t)$
\end{center}
Let $s\in \bigcup_{r\cdot t=p}(V_r\cap W_t)$. Then there are $r^*, t^*\in M$ such that $r^*\cdot t^*=p$ and $s\in V_{r^*}\cap W_{t^*}$. Since $s\in V_{r^*}$,\\
(1) There is an $A_1\subseteq M$ such that $A_1\ne \emptyset$, $supA_1=r^*$ and for all $q\in A_1$, $s \in U_q$.\\
(2) For all $B_1\subseteq M$ with $supB_1>r^*$, there is a $q\in B_1$ such that $s\notin U_q$.\\
Let $A_2=\{r\cdot t^*:r\in A_1\}$. Then $supA_2=\sum_{q\in A_1}q\cdot t^*=t^*\sum_{q\in A_1}q=t^*supA_1
=t^*\cdot r^*=p$. Let $q\in A_2$. There is a $r\in A_1$ such that $q=r\cdot t^*$ and hence, by (1),  $s\in U_r\cap W_{t^*}$. Therefore, $s\in\bigcup_{supA=p}\bigcap_{q\in A}\bigcup_{r\cdot t=q}(U_r\cap W_t)$.\\

Now let $s\in \bigcup_{supB>p}\bigcap_{q\in B}\bigcup_{r\cdot t=q}(U_r\cap W_t)$. (We show that this leads to a contradiction.) Therefore, there is a $B_1\subseteq M$ with $supB_1>p$ and\\	
(3) for all $q\in B_1$ there are $r_q,t_q\in M$ such that $r_q\cdot t_q=q$ and $s\in U_{r_q}\cap W_{t_q}$.\\
For all $q\in B_1$, $s\in W_{t_{q}}$. Since, $s\in W_{t_{q}}\cap W_{t^*}$, $s\in W_{t_{q}}\cap W_{t^*}\ne \emptyset$ and it follows that $t_q=t^*$. So we have\\
(4) for all $q\in B_1$, $t_q=t^*$ and hence, $r_q\cdot t^*=q$ and by (3), $s\in U_{r_{q}}$.\\
Let $C=\{r_q:q\in B_1\}$. Then\\
(5) $t^*\cdot supC=t^*\cdot\sum_{q\in B_1}r_q=\sum_{q\in B_1}t^*\cdot r_q=\sum_{q\in B_1}q=supB>p$.\\
Note that
\begin{center}
$sup(A_1\cup C)=\sum_{q\in A_1}q+\sum_{q\in B_1}r_q=supA_1+supC=r^*+supC\ge r^*$.
\end{center}
If $sup(A_1\cup C)>r^*$, then by (2), there is a $w\in A_1\cup C$ such that $s\notin U_w$.  By (1), $w\notin A_1$. Therefore, $w\in C$ and hence, $w=r_q$ for some $q\in B_1$ and we have $s\notin U_{r_q}$ which contradicts (4). So we have $sup(A_1\cup C)=r^*$. Therefore, $r^*+supC=supA_1+supC=Sup(A_1\cup C)=r^*$ and hence, $supC\le r^*$ and $t^*\cdot supC \le t^*\cdot r^*=p$. By (5) $t^*\cdot supC>p$ and hence, $p< p$, a contradiction. Therefore, $s\notin \bigcup_{supB>p}\bigcap_{q\in B}\bigcup_{r\cdot t=q}(U_r\cap W_t)$ and we have $s\in \bigcup_{supA=p}\bigcap_{q\in A}\bigcup_{r\cdot t=q}(U_r\cap W_t)\cdot-\bigcup_{supB>p}\bigcap_{q\in B}\bigcup_{r\cdot t=q}(U_r\cap W_t)$.\\
\\
Now let $s\in \bigcup_{supA=p}\bigcap_{q\in A}\bigcup_{r\cdot t=q}(U_r\cap W_t)\cdot-\bigcup_{supB>p}\bigcap_{q\in B}\bigcup_{r\cdot t=q}(U_r\cap W_t)$. Therefore,\\
(6) There is an $A_1\subseteq M$, $supA_1=p$ and for all $q\in A_1$, there are $r_q,t_q\in M$ such that $r_q\cdot t_q=q$ and $s\in U_{r_q}\cap W_{t_q}$\\
and\\
(7) for all $B_1\subseteq M$ with $sup B_1>p$, there is a $q\in B_1$ such that for all $r,t \in M$ with $r\cdot t=q$, $s\notin U_r\cap W_t$.\\
Since $s\in W_{t_q}$ for all $q\in A_1$, and $W_a\cap W_b=\emptyset$ for $a\ne b$, we have $t_{q_1}=t_{q_2}$ for all $q_1,q_2\in A_1$. Let $\bar{t}$ be the common value; that is, $t_q=\bar{t}$ for all $q\in A$. We then have\\
(8) for all $q\in A_1$, $s\in U_{r_q}$, $r_q\cdot \bar{t}=q$ and $s\in W_{\bar{t}}$.\\
We show that the assumption that for all $r\in M$ with $r\cdot \bar{t}=p$, $s\notin V_r$ leads to a contradiction and hence, $s\in V_{\bar{r}}$ for some $\bar{r}\in M$ such that $\bar{r}\cdot \bar{t}=p$ and $s\in \bigcup_{r\cdot t=p}(V_r\cap W_t)$.\\
Assume that for all $r\in M$ with $r\cdot\bar{t}=p$, $s\notin V_r$. Under this assumption we prove the following two lemmas.

\begin{Lemma}
\label{Lemma 1}
There is an $r_0 \in M$ such that $r_0\cdot \bar{t}=p$ and $s\in \bigcup_{supA=r_0}\bigcap_{q\in A}U_q$.
\begin{proof}
Let $A_2=\{r_q:q\in A_1\}$ and $r_0=supA_2$. Then\\
\begin{center}
$r_0\cdot \bar{t}=\bar{t}\cdot sup\{r_q:q\in A_1\}=\bar{t}\cdot \sum_{q\in A_1}r_q=\sum_{q\in A_1}\bar{t}\cdot r_q=\sum_{q\in A_1}q=supA_1=p$.
\end{center}
Let $w\in A_2$. Therefore, $w=r_q$ for some $q\in A_1$ and by (8), $s\in U_{r_q}$. Hence, $s\in \bigcup_{supA=r_0}\bigcap_{q\in A}U_q$.
\end{proof}
\end{Lemma}

\begin{Lemma}
\label{Lemma 2}
For all $r \in M$ with $r\cdot \bar{t}=p$, if $s\in\bigcup_{supA=r}\bigcap_{q\in A}U_q$, then there is an $r'\in M$ such that $r'>r$, $r'\cdot\bar{t}=p$ and $s\in\bigcup_{supA=r'}\bigcap_{q\in A}U_q$.
\begin{proof}
Let $r\in M$ such that $r\cdot\bar{t}=p$ and $s\in \bigcup_{supA=r}\bigcap_{q\in A}U_q$. We have assumed above that $s\notin V_r=\bigcup_{supA=r}\bigcap_{q\in A}U_q\cap-\bigcup_{supB>r}\bigcap_{q\in B}U_q$. Since $s\in \bigcup_{supA=r}\bigcap_{q\in A}U_q$, we have $s\in \bigcup_{supB>r}\bigcap_{q\in B}U_q$. Therefore,\\
(9) there is a $B_2\subseteq M$ such that $supB_2>r$ and for all $q\in B_2$, $s\in U_q$.\\Let $r'=supB_2$. 
Therefore, $r'>r$ and $\bar{t}\cdot supB_2\ge \bar{t} \cdot r=p$. We now show that we cannot have $\bar{t}\cdot supB_2>p$. Assume that $\bar{t}\cdot supB_2>p$. Let $C=\{\bar{t}\cdot q:q\in B_2\}$. Then $supC=\sum_{q\in B_2}\bar{t}\cdot q=\bar{t}\cdot \sum_{q\in B_2}q=\bar{t}\cdot supB_2>p$. By (7) there is a $z\in C$ such that for all $a,b\in M$ with $a\cdot b=z$, ,$s\notin U_a\cap W_b$. We have $z=\bar{t}\cdot q$ for some $q\in B_2$ and therefore, $s\notin U_q\cap W_{\bar{t}}$. However, $s\in W_{\bar{t}}$ by (8) and hence, $s\notin U_q$ which contradicts (9). We therefore have $\bar{t}\cdot r'=\bar{t}\cdot supB_2=p$, $r'>r$, $supB_2=r'$ and for all $q\in B_2$, $s\in U_q$. Hence, $s\in\bigcup_{supB=r'}\bigcap_{q\in B}U_q$.
\end{proof}
\end{Lemma}

We can now complete the proof of the theorem. Let $r_0$ be as in Lemma \ref{Lemma 1}. Suppose that $r_k$ has been defined and let $r_{k+1}=(r_k)'$ as in Lemma \ref{Lemma 2} . then $r_{k+1}>r_k$ for all $k\in \omega$ and since $M$ is finite this is impossible and this is our contradiction that proves that we cannot have for all $r\in M$ with $r\cdot\bar{t}=p$, $s\notin V_r$. Therefore, there is an $r\in M$ with $r\cdot \bar{t}=p$ such that $s\in V_r$; that is, $s\in\bigcup_{r\cdot t=p}(V_r\cap W_t)$.
\end{proof}

\end{Theorem}

\begin{Theorem}
\label{boolean 2}
Let $\mathfrak{B}=\langle B,+,\cdot ,-,0,1\rangle$ be a Boolean algebra, $p\in M$, and \\$Y:M\to B$. Then for all $p\in M$,\\
$\sum_{q\le p}(\prod_{supA=q}\prod_{r\in A}Y_r\cdot -\sum_{supA>q}\prod_{r\in A}Y_r)=(\sum_{q\le p}Y_q)\cdot (\prod_{q\nleq p}-Y_q)$\\
\begin{proof}
As before, let $Z$ be a subset of the set of subsets of a set $X$ such that $Z$ is closed under $\cup$, $\cap$, $-$ (complementation with respect to $X$), $\emptyset,X \in Z$, $p\in M$, and  $Y:M\to Z$. We show that\\
$\bigcup_{q\le p}(\bigcup_{supA=q}\bigcap_{r\in A}Y_r\cap -\bigcup_{supA>q}\bigcap_{r\in A}Y_r)=(\bigcup_{q\le p}Y_q)\cap (\bigcap_{q\nleq p}-Y_q)$\\
Let $s\in (\bigcup_{q\le p}Y_q)\cap (\bigcap_{q\nleq p}-Y_q)$ . There is a $q_0 \le p$ such that $s\in Y_{q_0}$ and for all $q\nleq p$, $s\notin Y_q$. Let
\begin{center}
$A_0=\{r:r\le p, s\in Y_r\}$
\end{center}
Note that $q_0\in A_0$ and hence, $A_0\ne \emptyset$. Let $q_1=supA_0$. For all $r\in A_0$, $r\le p$ and we have 

\begin{center}
$q_1= supA_0\le p$
\end{center}

Since $s\in Y_r$ for all $r\in A_0$, $s\in \bigcap_{r\in A_0}Y_r$. Therefore, 
\begin{center}
$s\in\bigcup_{supA=q_1}\bigcap_{r\in A}Y_r$.
\end{center}
Let $supA>q_1$.Then $A \nsubseteq A_0$. If so, $supA\le sup A_0=q_1$. Therefore, there is a $t\in A$ such that $t\notin A_0$.\\
If $t\le p$, then since $t\notin A_0$, $s\notin Y_t$.\\
If $t\nleq p$, then $s\notin Y_t$.\\
Therefore,  $s\notin Y_t$ and we have for all $A$ with $supA>q_1$, there is a $t\in A$ such that $s\notin Y_t$. Therefore, $s\in -\bigcup_{supA>q_1}\bigcap_{r\in A}Y_r$ and hence, $s\in \bigcup_{q\le p}(\bigcup_{supA=q}\bigcap_{r\in A}Y_r\cap -\bigcup_{supA>q}\bigcap_{r\in A}Y_r)$.\\

Let $s\in \bigcup_{q\le p}(\bigcup_{supA=q}\bigcap_{r\in A}Y_r\cap -\bigcup_{supA>q}\bigcap_{r\in A}Y_r)$.
There is a $q_0$ and there is an $A_0$ such that\\
1. $q_0\le p$\\
2. $A_0\ne \emptyset$\\
3. $supA_0=q_0$\\
4. for all $r\in A_0, s\in Y_r$\\
5. for all $A$ with $supA>q_0$, there is an $r\in A$ such that $s\notin Y_r$\\

Let $q\in A_0$. $q\le sup A_0=q_0\le p$. By 4, $s\in Y_q$ and hence, $s\in \bigcup_{q\le p}Y_q$.\\

Let $q \nleq p$.\\
Case (1) $q>p$. Then $sup\{q\}=q>p\ge q_0$. Therefore, $sup\{q\}>q_0$. By 5, $s\notin Y_q$.\\
Case (2) $q$ incomparable with $p$. We can't have $q\le q_0$; if so, then $q\le q_0\le p$. But $q$ is incomparable with $p$. Therefore, $q>q_0$ or $q$ incomparable with $q_0$.\\
If $q>q_0$, then $sup\{q\}=q>q_0$ and by 5, $s\notin Y_q$.\\
If $q$ is incomparable with $q_0$, let $B=A_0\cup\{q\}$. Then
\begin{center}
$supB=supA_0+q=q_0+q\ge q_0$.
\end{center}
If $q_0+q=q_0$, then $q\le q_0$. But $q$ is incomparable with $q_0$ and hence, $q_0+q>q_0$. Therefore, $supB>q_0$. By 5, there is an $r_0\in B$ such that $s\notin Y_{r_0}$. But $s\in Y_r$ for all $r\in A_0$ and hence $r_0\notin A_0$. Therefore, $r_0=q$ and $s\notin Y_q$.\\
In all cases, $s\notin Y_q$. and we have $s\in \bigcap_{q\nleq p}$. Therefore, $s\in (\bigcup_{q\le p}Y_q)\cap (\bigcap_{q\nleq p}-Y_q)$.
\end{proof}
\end{Theorem}

\begin{Theorem}
\label{is CA alpha}
If $\mathfrak{B}$ is a $\textrm{CA}_{\alpha}$, then $\mathfrak{M(B)}$ is an $\mathfrak{M}-\textrm{CA}_{\alpha}$.
\begin{proof}
Verifying Axioms 1 through 3 is straight forward. Using Theorem \ref{insert}, verifying Axioms 4, 5, 6, 10, 11, 12, and 13 through 31 is also straight forward. The verifications of Axioms 7, 8, and 9 follow.\\

Verification of Axiom 7: As before, define $Z_\kappa(A,x)=\prod_{p\in A}c_\kappa x^p$, where $x\in \mathfrak{|M(B)|}$. Then 
\begin{eqnarray*}
(\boldsymbol{c}_\kappa x)^r&=&\sum_{supA=r}Z_\kappa(A,x)\cdot -\sum_{supA>r}Z_\kappa(A,x)\\
&=&\sum_{supA=r}(Z_\kappa(A,x)\cdot -\prod_{supB>r}-Z_\kappa(B,x))
\end{eqnarray*}
Let $supA=r$ and $q\nleq r$. then $sup(A\cup \{q\})\geq r$. We can't have $sup(A\cup \{q\})= r$ since $q\nleq r$. Therefore,  $sup(A\cup \{q\})> r$.\\

\begin{eqnarray*}
Z_\kappa(A,x)\cdot \prod_{supB>r}-(Z_\kappa(B,x)&\le&Z_\kappa(A,x)\cdot -Z_\kappa(A\cup \{q\},x)\\
&=&Z_\kappa(A,x)\cdot -(Z_\kappa(A,x)\cdot Z_\kappa(\{q\},x))\\
&=&Z_\kappa(A,x)\cdot (-Z_\kappa(A,x)+ -Z_\kappa(\{q\},x))\\
&=&Z_\kappa(A,x)\cdot -Z_\kappa(A,x)+Z_\kappa(A,x)\cdot -Z_\kappa(\{q\},x)\\
&=&Z_\kappa(A,x)\cdot -Z_\kappa(\{q\},x)\\
&\le&-Z_\kappa(\{q\},x)\\
&=&-c_\kappa x^q\\
&\le&-x^q
\end{eqnarray*}
Therefore, $Z_\kappa(A,x)\cdot \prod_{SupB>r}-(Z_\kappa(B,x)\le -x^q$ for all $A$ with $supA=r$ and $q\nleq r$ and hence, $(\boldsymbol{c}_\kappa x)^r=\sum_{supA=p}(Z_\kappa(A,x)\cdot \prod_{supB>p}-Z_\kappa(B,x))\le -x^q$ for all $q\nleq r$. Therefore, $(\boldsymbol{c}_\kappa x)^r\le\prod_{q\nleq r}-x^q=-\sum_{q\nleq r}x^q$. Since $\sum_{q\nleq r}x^q+\sum_{q\le r}x^q=\sum_{q\in M}x^q=1$ and $\sum_{q\nleq r}x^q\cdot \sum_{q\le r}x^q=0$ we have $(\boldsymbol{c}_\kappa x)^r\le\sum_{q\le r}x^q$ and by Theorem \ref{less than}, $x\boldsymbol{+}\boldsymbol{c}_\kappa x=\boldsymbol{c}_\kappa x$.\\

Verification of Axiom 8: In a $\textrm{CA}_{\alpha}$, a $\kappa$-cylinder is a member a $x$ such that $c_\kappa x=x$. For the verification of Axiom 8a, we need the following:\\
If $x$ and $y$ are $\kappa$-cylinders, then $x+y$, $x\cdot y$, $-x$, and $c_\kappa x$ are $\kappa$-cylinders. It follows that if $\mathfrak{B}$ is a $\textrm{CA}_\alpha$ and $a\in \vert\mathfrak{M(B)}\vert$, then $(\boldsymbol{c}_\kappa a)^p=\sum_{supA=p}\prod_{q\in A}c_\kappa a^q\cdot -\sum_{supA>p}\prod_{q\in A}c_\kappa a^q$ is a $\kappa$-cylinder. By definition of $\boldsymbol{\cdot}$ in $\mathfrak{M(B)}$, 
$(\boldsymbol{c}_\kappa x\boldsymbol{\cdot}\boldsymbol{c}_\kappa y)^p=\sum_{q\cdot r=p}(\boldsymbol{c}_\kappa x)^q\cdot (\boldsymbol{c}_\kappa y)^r$ and

\begin{eqnarray*}
(\boldsymbol{c}_\kappa(x\boldsymbol{\cdot}\boldsymbol{c}_\kappa y))^p&=&\sum_{supA=p}\prod_{q\in A}c_\kappa(x\boldsymbol{\cdot}\boldsymbol{c}_\kappa y)^q\cdot -\sum_{supA>p}\prod_{q\in A}c_\kappa(x\boldsymbol{\cdot}\boldsymbol{c}_\kappa y)^q\\
&=&\sum_{supA=p}\prod_{q\in A}c_\kappa(\sum_{r\cdot t=q}x^r\cdot(\boldsymbol{c}_\kappa y)^t)\cdot -\sum_{supA=p}\prod_{q\in A}c_\kappa(\sum_{r\cdot t=q}x^r\cdot(\boldsymbol{c}_\kappa y)^t)\\
&=&\sum_{supA=p}\prod_{q\in A}\sum_{r\cdot t=q}c_\kappa (x^r\cdot(\boldsymbol{c}_\kappa y)^t)\cdot -\sum_{supA=p}\prod_{q\in A}\sum_{r\cdot t=q}c_\kappa (x^r\cdot(\boldsymbol{c}_\kappa y)^t)\\
&=&\sum_{supA=p}\prod_{q\in A}\sum_{r\cdot t=q} (c_\kappa x^r\cdot c_\kappa(\boldsymbol{c}_\kappa y)^t)\cdot -\sum_{supA=p}\prod_{q\in A}\sum_{r\cdot t=q} (c_\kappa x^r\cdot c_\kappa(\boldsymbol{c}_\kappa y)^t)\\
&=&\sum_{supA=p}\prod_{q\in A}\sum_{r\cdot t=q} (c_\kappa x^r\cdot(\boldsymbol{c}_\kappa y)^t)\cdot -\sum_{supA=p}\prod_{q\in A}\sum_{r\cdot t=q} (c_\kappa x^r\cdot(\boldsymbol{c}_\kappa y)^t)
\end{eqnarray*}
Since $(\boldsymbol{c}_\kappa x)^r=\sum_{supA=r}\prod_{a\in A}c_\kappa x^q\cdot-\sum_{supA>r}\prod_{a\in A}c_\kappa x^q$, using Theorem \ref{boolean 2}, where $U_q=c_\kappa x^q$, $V_r=(c_\kappa x)^r$, and $W_t=(c_\kappa y)^t$, we obtain $\boldsymbol{c}_\kappa (x\boldsymbol{\cdot} \boldsymbol{c}_\kappa y))^p=(\boldsymbol{c}_\kappa x\boldsymbol{\cdot} \boldsymbol{c}_\kappa y))^p$.\\

Verification of Axiom 9:\\ 

First we show, using Theorem  \ref{boolean 2}, that $\sum_{q\leq p}(\boldsymbol{c}_\kappa x)^q=c^\partial_\kappa\sum_{q\leq p}x^q$.\\
\begin{eqnarray*}
\sum_{q\leq p}(\boldsymbol{c}_\kappa x)^q&=&\sum_{q\leq p}(\sum_{supA=q}\prod_{r\in A}c_\kappa x^r\cdot -\sum_{supA>q}\prod_{r\in A}c_\kappa x^r)\\
&=&\sum_{q\leq p}c_\kappa x^q\cdot \prod_{q\nleq p}-c_\kappa x^q\\
&=&\sum_{q\leq p}c_\kappa x^q\cdot -\sum_{q\nleq p}c_\kappa x^q\\
&=&c_\kappa \sum_{q\leq p}x^q\cdot -c_\kappa \sum_{q\nleq p}x^q\\
&=&c_\kappa \sum_{q\leq p}x^q\cdot -c_\kappa -\sum_{q\leq p}x^q\\
&=&c_\kappa \sum_{q\leq p}x^q\cdot c^\partial_\kappa \sum_{q\leq p}x^q\\
&=&c^\partial_\kappa \sum_{q\leq p}x^q  \textrm{ since }  c^\partial_\kappa a\le c_\kappa a\\
\end{eqnarray*}
We then obtain $\sum_{q\le p}(\boldsymbol{c}_\kappa \boldsymbol{c}_\lambda x)^q=c^\partial_\kappa\sum_{q\le p}(\boldsymbol{c}_\lambda x)^q=c^\partial_\kappa c^\partial_\lambda \sum_{a\le p}x^q$. Similarily, $\sum_{q\le p}(\boldsymbol{c}_\lambda \boldsymbol{c}_\kappa x)^q=c^\partial_\lambda c^\partial_\kappa \sum_{a\le p}x^q$.\\

Therefore, $(\boldsymbol{c}_\kappa \boldsymbol{c}_\lambda x)^p\le \sum_{q\le p}(\boldsymbol{c}_\kappa \boldsymbol{c}_\lambda x)^q=c^\partial_\kappa c^\partial_\lambda\sum_{q\le p}x^q=c^\partial_\lambda c^\partial_\kappa\sum_{q\le p}x^q=\sum_{q\le p}(\boldsymbol{c}_\lambda \boldsymbol{c}_\kappa x)^p$. By Theorem \ref{less than}, $\boldsymbol{c}_ \lambda \boldsymbol{c}_ \kappa x\boldsymbol{+} \boldsymbol{c}_ \kappa \boldsymbol{c}_ \lambda x=\boldsymbol{c}_ \kappa \boldsymbol{c}_ \lambda x$. Similarly,  $\boldsymbol{c}_ \kappa \boldsymbol{c}_ \lambda  x\boldsymbol{+} \boldsymbol{c}_ \lambda \boldsymbol{c}_ \kappa x=\boldsymbol{c}_ \lambda \boldsymbol{c}_ \kappa x$. Using Axiom 1 (a), which is true in $\mathfrak{M(B)}$, we obtain $\boldsymbol{c}_\kappa \boldsymbol{c}_\lambda  x=\boldsymbol{c}_\lambda \boldsymbol{c}_\kappa  x$.
\end{proof}
\end{Theorem}

\begin{Definition}
Let $\mathfrak{A}=\langle A,+,\cdot ,-,u_p,c_\kappa ,d_{\kappa,\lambda},\delta_p\rangle_{\kappa,\lambda<\alpha,p\in M}$ be an $\mathfrak{M}-\textrm{CA}_{\alpha}$. $\mathfrak{C(A)}=\langle C, +, \cdot, -, U_0, U_1, c_\kappa, d_{\kappa,\lambda}\rangle_{\kappa,\lambda<\alpha}$, where $C=\{\delta_1x:x\in \vert \mathfrak{A}\vert\}$ and all operations in $\mathfrak{C(A)}$ are the restrictions of the operations in $\mathfrak{A}$. Note that by Theorem \ref{insert}, $C$ is closed under all operations.

\end{Definition}

\begin{Theorem}
If $\mathfrak{A}$ an $\mathfrak{M}-\textrm{CA}_{\alpha}$, then $\mathfrak{C(A)}$ is a $\textrm{CA}_{\alpha}$.

\begin{proof}
The proof follows easily from the axioms for $\mathfrak{M}-\textrm{CA}_{\alpha}$.
\end{proof}
\end{Theorem}

\begin{Definition}
If $\mathfrak{A}$ and $\mathfrak{B }$ are algebras then $\mathfrak{A}\preccurlyeq \mathfrak{B}$ if $\mathfrak{A}$ is a subalgebra of $\mathfrak{B}$ and $\mathfrak{A}\precsim \mathfrak{B}$ if $\mathfrak{A}$ is isomorphic to a subalgebra of $\mathfrak{B}$.
\end{Definition}

\begin{Theorem}
\label{embed}
If $\mathfrak{A}$ is an $\mathfrak{M}-\textrm{CA}_{\alpha}$, then $\mathfrak{A}\precsim \mathfrak{M(C(A))}$.

\begin{proof}
Let $\mathfrak{A}$ be an $\mathfrak{M}-\textrm{CA}_{\alpha}$. Define $f:\vert \mathfrak{A}\vert \to \vert \mathfrak{M(C(A))} \vert$ by $f(a)^p=\delta_p(a)$ for $a\in \vert\mathfrak{A}\vert$. It is easily verified that $f$ is an isomorphism from  $\mathfrak{A}$ into $\mathfrak{M(C(A))}$.

\end{proof}
\end{Theorem}

If $\mathfrak{C}_j$ are algebras for $j\in J$, then $\prod\langle\mathfrak{C}_j:j\in J\rangle$ is the product of the algebras $\mathfrak{C}_j$  for $j\in J$

\begin{Theorem}
\label{prod 1}
If $\mathfrak{B}$ and $\mathfrak{C}_j$ for $j\in J$ are $\textrm{CA}_{\alpha}$ and $\mathfrak{B}\precsim \prod\langle\mathfrak{C}_j:j\in J\rangle$, then $\mathfrak{M(B)}\precsim \prod\langle\mathfrak{M(C_j)}:j\in J\rangle$.
\begin{proof}
Let $f:\vert\mathfrak{B}\vert\to \prod\langle\vert\mathfrak{C}_j\vert:j\in J\rangle$ be a monomorphism. Define $g:\vert\mathfrak{M(B)}\vert\to \prod\langle\vert\mathfrak{M(C_j)}\vert:j\in J\rangle$ by $(g(a)_j)^p=f(a^p)_j$ for $a\in \vert\mathfrak{M(B)}\vert$. It is tedious but straight forward to verify that $g$ is a monomorphism from $\mathfrak{M(B)}$ into $\prod\langle\mathfrak{M(C_j)}:j\in J\rangle$.
\end{proof}

\end{Theorem}

\begin{Theorem}
If $\mathfrak{B}$ is a $\textrm{CA}_{\alpha}$, then $\mathfrak{C(M(B))}\simeq\mathfrak{B}$.
\label{iso}
\begin{proof}
Define $g:\vert\mathfrak{B}\vert\to \vert\mathfrak{M(C)}\vert$ by
\begin{displaymath}
g(a)^p=\left \{ \begin{array}{ll}
a & \textrm{if $p=1$}\\
0 &  \textrm{if $p\ne 0,1$}\\
-a & \textrm{if $p=0$}\end{array} \right .
\end{displaymath}
Define $f:\vert\mathfrak{B}\vert\to \vert\mathfrak{C(M(B))}\vert$ by $f(a)=\delta_1g(a)$. Using Theorem \ref{insert}, it it easily shown that $f$ is an isomorphism from $\mathfrak{B}$ to $\mathfrak{C(M(B))}$.
\end{proof}
\end{Theorem}

\begin{Theorem}
\label{prod 2}
If $\mathfrak{A}$ and $\mathfrak{D}_j$ are $\mathfrak{M}-\textrm{CA}_{\alpha}$ for $j\in J$ such that $\mathfrak{A}\precsim\prod\langle\mathfrak{D}_j:j\in J\rangle$, then $\mathfrak{C(A)}\precsim\prod\langle\mathfrak{C(D})_j:j\in J \rangle$.

\begin{proof}
Let $f:\vert\mathfrak{A}\vert\to \prod\langle\vert\mathfrak{D}_j\vert:j\in J\rangle$ be a monomorphism. Let $f=g\lceil \vert\mathfrak{C(A)}\vert$, the restriction of $g$ to $\vert\mathfrak{C(A)}\vert$. It is easily verified that $g$ is a monomorphism from $\mathfrak{C(A)}$ to $\prod\langle\mathfrak{C(D}_j):j\in J\rangle$.
\end{proof}

\end{Theorem}

\begin{Theorem}
\label{C embed}
Let $\mathfrak{A}$ and $\mathfrak{B}$ be $\mathfrak{M}-\textrm{CA}_{\alpha}$ such that $\mathfrak{A}\preccurlyeq\mathfrak{B}$. Then $\mathfrak{C(A)}\preccurlyeq\mathfrak{C(B)}$. 
\begin{proof}
The proof follows easily from the definition of $\mathfrak{C(A)}$.
\end{proof}
\end{Theorem}

\subsection{$\mathfrak{M}$-cylindric set algebras}

\begin{Definition}
$\mathfrak{A}$ is a $\mathfrak{M}$-cylindric set algebra of dimension $\alpha$ with base $U$ ($\mathfrak{M}-\textrm{CSA}_{\alpha}$) if $\mathfrak{A}\preccurlyeq\mathfrak{M(B)}$ for some cylindric set algebra $\mathfrak{B}$ of dimension $\alpha$ with base $U$. $\mathfrak{A}$ is full if, in addition, $|\mathfrak{B}|=\mathcal{P}(U^\alpha)$.
\end{Definition}

We will use the following notation to indicate that we are working in a $\mathfrak{M}-\textrm{CSA}_{\alpha}$:

$\sqcup$ for $\boldsymbol{+}$, $\sqcap$ for $\boldsymbol{\cdot}$, $\circleddash$ for $\boldsymbol{-}$, $\boldsymbol{C}_\kappa$ for $\boldsymbol{c}_\kappa$, $\boldsymbol{D}_{\kappa, \lambda}$ for $\boldsymbol{d}_{\kappa, \lambda}$, $\boldsymbol{U}_p$ for $\boldsymbol{u}_p$, $\boldsymbol{\Delta}_p$ for $\boldsymbol{\delta}_p$ and $\bigsqcup$ for $\sum$.\\

In addition, if $l$ is a member of a $\mathfrak{M}-\textrm{CSA}_{\alpha}$, then $\boldsymbol{\Delta}L=\boldsymbol{\delta}_0L \sqcup \boldsymbol{\delta}_1L$.

\begin{Theorem}
If $\mathfrak{A}$ is a $\mathfrak{M}-\textrm{CSA}_{\alpha}$, then $\mathfrak{C(A)}\simeq \mathfrak{B}$ for some cylindric set algebra $\mathfrak{B}$. 
\begin{proof}
Let $\mathfrak{A}$ be a $\mathfrak{M}-\textrm{CSA}_{\alpha}$. Then $\mathfrak{A}\preccurlyeq\mathfrak{M(D)}$ for some cylindric set algebra $\mathfrak{D}$. By Theorems \ref{C embed} and \ref{iso}, $\mathfrak{C(A)}\preccurlyeq\mathfrak{C(M(D))}\simeq\mathfrak{D}$. Therefore, $\mathfrak{C(A))}$ is isomorphic to a cylindric set algebra.
\end{proof}
\end{Theorem}

If $\mathfrak{A}$ is a $\mathfrak{M}-\textrm{CSA}_{\alpha}$, then $\mathfrak{A}\preccurlyeq\mathfrak{M(B)}$ for some cylindric set algebra $\mathfrak{B}$ of dimension $\alpha$. If the base of $\mathfrak{B}$ is $U$, that is, $\vert\mathfrak{B}\vert \subseteq\mathcal{P}(U^\alpha)$,  where $\mathcal{P}(U^\alpha)$ is the set of sunsets of $U^\alpha$, then $\vert\mathfrak{A}\vert\subseteq \{L:L\in\mathcal{P}(U^\alpha)^M  \textrm{ and }  \bigcup_{p\in M}L^p=U  \textrm{ and if }  p\ne q  \textrm{ then }  L^p\cap L^q=\emptyset\}$.\\

\begin{Definition} For $s\in U^\alpha$ and $x\in U$ $s\binom{\kappa}{x}\in U^\alpha$ is defined by
\begin{displaymath}
s\binom{\kappa}{x}_\lambda=\left \{ \begin{array}{ll}
s_\lambda & \textrm{if $\kappa\ne \lambda $}\\
x &  \textrm{if $\kappa=\lambda$}\end{array} \right .
\end{displaymath}
\end{Definition}

\begin{Definition} 
For $X\subseteq U^\alpha$, and $y\in U$, 
\begin{center}
$y\circ_\kappa X=\{s:s\in U^\alpha \textrm{ and } s\binom{\kappa}{y}\in X\}$
\end{center}
\end{Definition}

For $X\subseteq U^\alpha$, $C_\kappa X=\{s:s\in U^\alpha \textrm{ and for some }x\in U, s\binom{\kappa}{x}\in X\}$, the usual cylindrification  in cylindric set algebras.

For $L\in\mathcal{P}(U^\alpha)^M$, we say that $L$ partitions $U^\alpha$ if $\{L^p:p\in M\}$ is a partition of $U^\alpha$.

\begin{Theorem}
\label{motivation 1}
If $x\in U$ and $L\in\mathcal{P}(U^\alpha)^M$ partitions $U^\alpha$, then $\{x\circ_\kappa L^p:p\in M\}$ partitions $U^\alpha$.
\end{Theorem}

\begin{Definition} 
For $j\in M^U$, let $\hat{j}=sup \{j_x:x\in U\}$.
\end{Definition}

\begin{Definition} 
For $L\in\mathcal{P}(U^\alpha)^M$ where $L$ partitions $U^\alpha$ define $t_\kappa L\in \mathcal{P}(U^\alpha)^M$ by $(E_\kappa L)^p=\bigcup_{j\in M^U,\hat{j}=p}\bigcap_{x\in U} x\circ_\kappa L^{j_x}$.
\end{Definition}

We will show that $E_\kappa L=\boldsymbol{C}_\kappa L$, the cylindrification in $\mathfrak{M(B)}$ where $\mathfrak{B}$ is a $\textrm{CSA}_{\alpha}$.

\begin{Theorem}
\label{motivation 2}
If  $L\in\mathcal{P}(U^\alpha)^M$ partitions $U^\alpha$, $\kappa<\alpha$, $p\in M$, $s\in (\boldsymbol{C}_\kappa L)^p$ and $s\in \bigcap_{x\in U}x\circ_\kappa L^{j_x}$ where $j\in M^U$, then $\hat{j}=p$.
\begin{proof}
Since $s\in (\boldsymbol{C}_\kappa L)^p=\bigcup_{supA=p}\bigcap_{q\in A}C_\kappa L^q\cdot -\bigcup_{supA>p}\bigcap_{q\in A}C_\kappa L^q$, there is an $A\subseteq M$ with $supA=p$ such that for all $q\in A$, $s\in C_\kappa L^q$. Therefore, for all $q\in A$, there is an $x_q\in U$ such that $s\binom{\kappa}{x_q}\in L^q$; that is, $s\in x_q\circ_\kappa L^q$. Therefore, for all $q\in A$, $s\in x_q\circ_\kappa L^q\cap x_q\circ_\kappa L^{j_{x_q}}$. By Theorem \ref{motivation 1}, $j_{x_q}=q$ for all $q\in A$. Therefore, $\hat{j}=sup\{j_x:x\in U\}\ge sup\{j_{x_q}:q\in A\}=sup\{q:q\in A\}=supA=p$ and hence $\hat{j}\ge p$.\\

Suppose that $\hat{j}>p$. Let $B=\{j_x:x\in U\}$. Then $supB =\hat{j}>p$. Since $s\notin \bigcup_{supA>p}\bigcap_{q\in A}C_\kappa L^q$, we have $s\notin C_\kappa L^q$ for some $q\in B $ and hence, $s\binom{\kappa}{y}\notin L^q$ for all $y\in U$. Therefore, $s\notin y\circ_\kappa L^q$ for all $y\in U$. Since $q\in B$, $q=j_x$ for some $x\in U$. Therefore, $s\notin y\circ_\kappa L^{j_x}$ for all $y\in U$. In particular,  $s\notin x\circ_\kappa L^{j_x}$, and this contradicts $s\in\bigcap_{x\in U}x\circ_\kappa L^{j_x}$. Therefore, $\hat{j} =p$
\end{proof}
\end{Theorem}

\begin{Theorem}
\label{motivation 3}
If $L\in\mathcal{P}(U^\alpha)^M$ partitions $U^\alpha$, then $E_\kappa L$ partitions $U^\alpha$.
\end{Theorem}

\begin{Theorem}
\label{motivation 4}
If $L\in\mathcal{P}(U^\alpha)^M$ partitions $U^\alpha$ and $s\in U^\alpha$, then there is a unique $j\in M^U$ such that $s\in  \bigcap_{x\in U}x\circ_\kappa L^{j_x}$.
\begin{proof}
This follows easily from Theorem \ref{motivation 1}.
\end{proof}
\end{Theorem}

\begin{Theorem}
If $L\in\mathcal{P}(U^\alpha)^M$ partitions $U^\alpha$, then $E_\kappa L=\boldsymbol{C}_\kappa L$.
\label{motivation 5}
\begin{proof}
Assume that $L\in\mathcal{P}(U^\alpha)^M$ partitions $U^\alpha$.\\Let $s\in (E_\kappa L)^p=\bigcup_{j\in M^U,\hat{j}=p}\bigcap_{x\in U} x\circ_\kappa L^{j_x}$. Therefore, $s\in\bigcap_{x\in U}x\circ_\kappa L^{j_x}$ for some $j\in M^U$ with $\hat{j}=p$. Let $A=rng(j)$, Therefore, $supA=\hat{j}=p$. Let $q\in A$. Then $q=j_x$ for some $x\in U$. Therefore, $s\in x\circ_\kappa L^{j_x}=x\circ_\kappa L^q$; that is, $s\binom{\kappa}{x}\in L^q$. Therefore, $s\in C_\kappa L^q$. We have shown that for all $q\in A$, $s\in C_\kappa L^q$ and hence, $s\in\bigcap_{q\in A}C_\kappa L^q$ and it follows that $s\in\bigcup_{supA=p}\bigcap_{q\in A}C_\kappa L^q$.\\
Now suppose that $s\in\bigcup_{supA>p}\bigcap_{q\in A}C_\kappa L^q$. Then there is an $A\subseteq M$, with $supA>p$ and for all $q\in A$, $s\in C_\kappa L^q$. Therefore, for all $q\in A$, there is an $x_q\in U$ such that $s\in x_q\circ_\kappa L^q$. By Theorem \ref{motivation 4}, there is a unique $j\in M^U$ such that $s\in \bigcap_{x\in U}x\circ_\kappa L^{j_x}$. Therefore, $s\in x_q\circ_\kappa L^{j_{x_q}}$ for all $q\in A$ and it follows that for all $q\in A$, $j_{x_q}=q$. Let $t=\hat{j}$. Since $s\in \bigcap_{x\in U}x\circ_\kappa L^{j_x}$, $s\in \bigcup_{\hat{j}=t}\bigcap_{x\in U}x\circ_\kappa L^{j_x}=(E_\kappa L)^t$. Since $s\in (E_\kappa L)^p$, by Theorem \ref{motivation 3}, $t=p$. We then have\\
$p<supA=sup\{q:q\in A\}=sup\{j_{x_q}:q\in A\}\leq sup \textrm{ }rng(j)=\hat{j}=p$. Hence, $p<p$ which is a contradiction.\\
Therefore we have $s\notin\bigcup_{supA>p}\bigcap_{q\in A}C_\kappa L^q$ and hence,\\
 $s\in \bigcup_{supA=p}\bigcap_{q\in A}C_\kappa L^q\cdot -\bigcup_{supA>p}\bigcap_{q\in A}C_\kappa L^q=(\boldsymbol{C}_\kappa L)^p$.\\
 
 Let $s\in(\boldsymbol{C}_\kappa L)^p$. By Theorem \ref{motivation 4}, there is a $j\in M^U$ such that $s\in \bigcap_{x\in U}x\circ_\kappa L^{j_x}$. By Theorem \ref{motivation 2}, $\hat{j}=p$. Therefore, $s\in \bigcup_{\hat{j}=p}\bigcap_{x\in U}x\circ_\kappa L^{j_x}=(E_\kappa L)^p$.
\end{proof}
\end{Theorem}

\begin{Theorem} In an $\mathfrak{M}-\textrm{CSA}_{\alpha}$ with base $U$, $\boldsymbol{C}_\kappa (L\sqcup K)=\boldsymbol{C}_\kappa L\sqcup \boldsymbol{C}_\kappa K$.
\label{motivation 6}
\begin{proof}
We use Theorem \ref{motivation 5} and show that $E_\kappa (L\sqcup K)=E_\kappa L \sqcup E_\kappa K$.
Let $s\in (E_\kappa(L\sqcup K))^p$. Then there is a $j\in M^U$ with $\hat{j}=p$ and for all $x\in U$, $s\in x\circ_\kappa(L\sqcup K)^{j_x}=\bigcup_{r+t=j_x}L^r\cap K^t$. Therefore, $s\binom{\kappa}{x}\in (L\sqcup K)^{j_x}$ for all $x\in U$. There are $r_x$, $t_x \in M$, such that $r_x+t_x=j_x$ and $s\in L^{r_x}\cap K^{t_x}$ for all $x\in U$.\\
Let $m=\hat{r}$ and $n=\hat{t}$. Then $\hat{r}+\hat{t}=\sum_{x\in U}r_x+\sum_{x\in U}t_x=\sum_{x\in U}(r_x+t_x)=\sum_{x\in U}j_x=\hat{j}$. Therefore, $m+n=\hat{r}+\hat{t}=\hat{j}=p$. We have 
$\hat{r}=m$ and for all $x\in U$, $s\in x\circ_\kappa L^{r_x}$ and hence, $s\in \bigcup_{\hat{r}=m}\bigcap_{x\in U}x\circ_\kappa L^{r_x}=(\boldsymbol{c}_\kappa L)^m$. Similarly,  $s\in (\boldsymbol{c}_\kappa K)^n$. Therefore, $s\in \bigcup_{m+n=p}(E_\kappa L)^m\cap(E_\kappa K)^n=(E_\kappa L\cup E_\kappa K)^p$,\\

Now let $s\in (E_\kappa L \sqcup E_\kappa K)^p$. There are $r$, $t\in M$ with $r+t=p$ and $s\in (E_\kappa L)^r$ and $s\in (E_\kappa K)^t$. There are $h\in M^U$ and $i\in M^U$ with $\hat{h}=r$ and $\hat{i}=t$ such that $s\in x\circ_\kappa L^{h_x}$ and $s\in x\circ_\kappa K^{i_x}$ for all $x\in U$. Let $j_x=h_x+i_x$, Then $\hat{j}=\hat{h}+\hat{i}=r+t=p$. Let $x\in U$, $a=h_x$ and $b=i_x$. Then $a+b=h_x+i_x=j_x$ and $s\in x\circ_\kappa L^{h_x}=x\circ_\kappa L^a$ and $s\in x\circ_\kappa K^{i_x}=x\circ_\kappa K^b$. Therefore, $s\binom{\kappa}{x} \in L^a\cap K^b$ and hence,  $s\binom{\kappa}{x} \in \bigcup_{a+b=j_x}L^a\cap K^b=(L\sqcup K)^{j_x}$. Therefore, $s\in x\circ_\kappa (L\sqcup K)^{j_x}$. Since $\hat{j}=p$, $s\in (E_\kappa (L\sqcup K))^p$.\\
Therefore, $(E_\kappa(L\cup K))^p=(E_\kappa L \cup E_\kappa K)^p$ for all $p\in M$ and $E_\kappa(L\cup K)=E_\kappa L \cup E_\kappa K$.
\end{proof}.
\end{Theorem}

\begin{Theorem} In an $\mathfrak{M}-\textrm{CSA}_{\alpha}$, if $L\le K$ then $\boldsymbol{C}_\kappa L\le \boldsymbol{C}_\kappa K$.
\label{c is monotone}
\end{Theorem}

\subsection{Dimension sets}

Recall that in a $\textrm{CA}_{\alpha}$ $\mathfrak{A}$, $dim_{\mathfrak{A}}a=\{\kappa:c_\kappa a\ne a\}$.

\begin{Definition}
Let $a$ be an element in an $\mathfrak{M}-\textrm{CA}_{\alpha}$ $\mathfrak{A}$. $dim_{\mathfrak{A}}a=\bigcup_{p\in M}\{\kappa:c_\kappa\delta_pa\ne\delta_pa\}$.
\end{Definition}

\begin{Theorem}
\label{dim 1}
If $\mathfrak{B}$ is a $\textrm{CA}_{\alpha}$ and $a\in |\mathfrak{M(B)}|$, then $dim_{\mathfrak{M(B)}}a=\bigcup_{p\in M}dim_{\mathfrak{B}}a^p$.
\begin{proof}
Let $\kappa\notin dim_{\mathfrak{M(B)}}a$. Therefore, $\boldsymbol{c}_{\kappa}\boldsymbol{\delta}_pa=\boldsymbol{\delta}_pa$ for all $p \in M$. By definition of $\delta_{p}$ in $\mathfrak{M(B)}$ and Theorem \ref{insert}, $c_\kappa a^p=c_\kappa (\boldsymbol{\delta}_pa)^1=(\boldsymbol{c}_{\kappa}\boldsymbol{\delta}_pa)^1=(\boldsymbol{\delta}_pa)^1=a^p$ for all $p \in M$. Therefore, $\kappa\notin dim_{\mathfrak{B}}a^p$ for all $p \in M$.

Let $\kappa\notin dim_{\mathfrak{B}}a^p$ for all $p \in M$. Therefore, $c_\kappa a^p=a^p$ for all $p \in M$ and hence,  by Theorem \ref{insert}, $(\boldsymbol{c}_{\kappa}\boldsymbol{\delta}_pa)^1=c_\kappa (\boldsymbol{\delta}_pa)^1=c_\kappa a^p=a^p=(\boldsymbol{\delta}_pa)^1$.
\end{proof}
\end{Theorem}

\begin{Definition}
An $\mathfrak{M}-\textrm{CA}_{\alpha}$ $\mathfrak{A}$ is locally finite if $dim_{\mathfrak{A}}a$ is finite for all $a\in \vert\mathfrak{A}\vert$.
\end{Definition}

\begin{Theorem}
If $\mathfrak{B}$ is a locally finite $\textrm{CA}_{\alpha}$, then $\mathfrak{M(B)}$ is a locally finite $\mathfrak{M}-\textrm{CA}_{\alpha}$.
\begin{proof}
The proof follows easily from Theorem \ref{dim 1}.
\end{proof}
\end{Theorem}

\begin{Theorem}
\label{dim 2}
Let $\mathfrak{A}$ be a $\mathfrak{M}-\textrm{CA}_{\alpha}$ and $a\in \vert\mathfrak{A}\vert$. If $\kappa\notin dim_\mathfrak{A}a$, then $c_\kappa a=a$; that is, $\{\kappa:c_\kappa a \neq a\} \subseteq dim_\mathfrak{A}a$.
\begin{proof}
Assume that $\kappa \notin dim_\mathfrak{A}a$. Therefore, $c_\kappa \delta_p a=\delta_pa$ for all $p\in M$.\\
By Axiom 31 we have
\begin{eqnarray*}
\delta_p(c_\kappa a)&=&\sum_{sup A=p}\prod_{q\in A}c_\kappa\delta_qa\cdot-\sum_{sup A>p}\prod_{q\in A}c_\kappa\delta_qa\\
&=&\sum_{sup A=p}\prod_{q\in A}\delta_qa\cdot-\sum_{sup A>p}\prod_{q\in A}\delta_qa\\
&=&\delta_pa\cdot-\sum_{r>p}\delta_ra \textrm{   since } \delta_qa\cdot\delta_pa=0 \textrm{ if }q\neq p\\
&=&\delta_pa\cdot\prod_{r>p}-\delta_ra\\
&=&\delta_pa
\end{eqnarray*}
The last line follows from the fact that if $r>p$, then $r\ne p$ and hence, $\delta_pa\cdot\delta_ra=0$ and therefore, $\delta_pa\leq-\delta_ra$

Therefore, for all $p\in M$, $\delta_p(c_\kappa a)=\delta_pa$ and by Axiom 30, $c_\kappa a=a$.
\end{proof}
\end{Theorem}

\begin{Theorem}
\label{dim 3}
If $\mathfrak{A}$ is a locally finite $\mathfrak{M}-\textrm{CA}_{\alpha}$, then $\mathfrak{C(A)}$ is a locally finite $\textrm{CA}_{\alpha}$.
\begin{proof}
Note that $c^\mathfrak{C(A)}_\kappa a=c^\mathfrak{A}_\kappa a$. By Theorem \ref{dim 2}, $dim_\mathfrak{C(A)}a=\{\kappa:c^\mathfrak{C(A)}_\kappa a\neq a\}=\{\kappa:c^\mathfrak{A)}_\kappa a\neq a\}\subseteq dim_\mathfrak{A}a$.
\end{proof}
\end{Theorem}

\subsection{Representation of $\mathfrak{M}-\textrm{CA}_{\alpha}'s$}

\begin{Definition}
A member $A$ of a $\textrm{CSA}_{\alpha}$ with base $U$ depends on $\Gamma\subseteq \alpha$ if for all $s$, $t\in U^\alpha$, $s\lceil \Gamma=t\lceil \Gamma$ implies that $s\in A$ if and only if $t\in A$. A member $L$ of a $\mathfrak{M}-\textrm{CSA}_{\alpha}$ depends on $\Gamma\subseteq \alpha$ if for all $p\in M$, $L^P$ depends on $\Gamma$.
\end{Definition}

\begin{Definition}
A $\textrm{CSA}_{\alpha}$ is regular if every member depends on its dimension set. A $\mathfrak{M}-\textrm{CSA}_{\alpha}$ is regular if every member depends on its dimension set.
\end{Definition}

\begin{Theorem}
\label{regular}
If a $\textrm{CSA}_{\alpha}$ $\mathfrak{B}$ is regular, then $\mathfrak{M(B)}$ is regular.
\begin{proof}
The proof follows directly from Theorem \ref{dim 1}.
\end{proof}
\end{Theorem}

\begin{Definition}
A $\textrm{CA}_{\alpha}$  is representable if it is isomorphic to a subdirect product of $\textrm{CSA}_{\alpha}$'s.
A $\mathfrak{M}-\textrm{CA}_{\alpha}$ is representable if it is isomorphic to a subdirect product of $\mathfrak{M}-\textrm{CSA}_{\alpha}$'s.
\end{Definition}

A subalgebra of a $\mathfrak{M}-\textrm{CA}_{\alpha}$ is a $\mathfrak{M}-\textrm{CA}_{\alpha}$ since all the axioms are universal. Therefore, a $\mathfrak{M}-\textrm{CA}_{\alpha}$ is representable if it is isomorphic to a subalgebra of a direct product of $\mathfrak{M}-\textrm{CSA}_{\alpha}$'s.

\begin{Theorem}
\label{rep 1}
If $\mathfrak{A}$ is a $\mathfrak{M}-\textrm{CA}_{\alpha}$ and $\mathfrak{C(A)}$ is isomorphic to a subalgebra of a product of (regular) $\textrm{CSA}_{\alpha}$'s, then $\mathfrak{A}$ is isomorphic to a subalgebra of a product of (regular) $\mathfrak{M}-\textrm{CSA}_{\alpha}$'s and hence, if $\mathfrak{C(A)}$ is representable,then so is $\mathfrak{A}$.
\begin{proof}
Assume that $\mathfrak{C(A)}\precsim\prod\langle\mathfrak{C}_j:j\in J\rangle$ where $\mathfrak{C}_j$ is a (regular) $\textrm{CSA}_{\alpha}$ for $j\in J$. By Theorem \ref{prod 1}, $\mathfrak{M(C(A))}\precsim\prod\langle\mathfrak{M(C_j)}:j\in J\rangle$. (By Theorem \ref{regular} $\mathfrak{M(C_j)}$ is regular.) By Theorem \ref{embed},  $\mathfrak{A}\precsim\mathfrak{M(C(A))}$. Therefore, $\mathfrak{A}\precsim\prod\langle\mathfrak{M(C_j)}:j\in J\rangle$.
\end{proof}
\end{Theorem}

\begin{Theorem} (Representation Theorem)
\label{rep 2}
Every locally finite $\mathfrak{M}-\textrm{CA}_{\alpha}$ is isomorphic to a subalgebra of a product of regular $\mathfrak{M}-\textrm{CSA}_{\alpha}$ and hence, is representable.
\begin{proof}
Let $\mathfrak{A}$ be a locally finite $\mathfrak{M}-\textrm{CA}_{\alpha}$. By Theorem \ref{dim 3}, $\mathfrak{C(A)}$ is a locally finite $\textrm{CA}_{\alpha}$. Therefore, $\mathfrak{C(A)}$ is representable [4]; that is, $\mathfrak{C(A)}$ is isomorphic to a subalgebra of a product of regular $\textrm{CSA}_{\alpha}$'s. By Theorem \ref{rep 1}, $\mathfrak{A}$ s isomorphic to a subalgebra of a product of regular $\mathfrak{M}-\textrm{CSA}_{\alpha}$ and hence is representable.
\end{proof}
\end{Theorem}

\begin{Theorem}
\label{rep 3}
If $\mathfrak{A}$ is a representable $\mathfrak{M}-\textrm{CA}_{\alpha}$, then $\mathfrak{C(A)}$ is a representable $\textrm{CA}_{\alpha}$.
\begin{proof}
Let $\mathfrak{A}$ be a representable $\mathfrak{M}-\textrm{CA}_{\alpha}$. Therefore, $\mathfrak{A}\precsim\prod\langle\mathfrak{D_j}:j\in J\rangle$ where $\mathfrak{D_j}$ is a $\mathfrak{M}-\textrm{CSA}_{\alpha}$ for $j\in J$. By Theorem \ref{prod 2} $\mathfrak{C(A)}\precsim\prod\langle\mathfrak{C(D_j)}:j\in J\rangle$. $\mathfrak{D_j}=\mathfrak{M(B_j)}$ where $\mathfrak{B_j}$ is a $\textrm{CSA}_{\alpha}$. By Theorem \ref{iso}, $\mathfrak{C(D_j)}= \mathfrak{C(M(B_j))}\simeq\mathfrak{B_j}$. Therefore, for all $j\in J$, $\mathfrak{C(D_j)}$ is isomorphic to a $\textrm{CSA}_{\alpha}$. Therefore, $\mathfrak{C(A)}$ is representable.
\end{proof}
\end{Theorem}

Let $\mathfrak{B}$ be a non-representable $\textrm{CA}_{\alpha}$ [3]. If $\mathfrak{M(B))}$ is representable then by Theorem \ref{rep 3}, $\mathfrak{C(M(B))}$ is representable. By Theorem \ref{iso}, $\mathfrak{C(M(B))}\simeq \mathfrak{B}$. Therefore, $\mathfrak{B}$ is representable, which is a contradiction. Therefore, there is a non-representable $\mathfrak{M}-\textrm{CA}_{\alpha}$.

\section{$\mathfrak{M}$-Valued Logic}

\subsection{Syntax of $\mathfrak{M}$-Valued Logic}

The symbols of this language are variables, $v_k$ for $k\in\omega$, the set of natural numbers; $n_i$-ary relation symbols $R_i$ for $i\in I$ where $I$ is a nonempty set; equal symbol $\thickapprox$; negation symbol $\neg$; and symbol $\wedge$; or symbol $\vee$; existential qunatifier $\exists$; and two additional symbols, $\gamma_p$ and  $t_p$ for $p\in M$.
Formulas are defined inductively as follows:
\begin{Definition}

$R_i v_{k_0} \cdots v_{k_{n_i-1}}$, $v_j \thickapprox v_k$, and $t_p$  are formulas for $i\in I,j,k\in \omega$, and $p\in M$.

If $\phi$ and $\theta$ are formulas, then $\neg \phi$, $( \phi \wedge \theta)$, $ (\phi \vee \theta)$, $\exists v_k \phi$, and  $\gamma_p \phi$ are all formulas.
\end{Definition}

\begin{Definition} 
The set of free variables of a formula $\phi$, $Fv(\phi)$,  is defined inductively as follows:
\begin{tabbing}
xxxxxxxxxxxxxxxxxxxxxxxxxxxxxxxxxxxx\= \kill

$Fv(t_p)=\emptyset$ \> $Fv(v_i\thickapprox v_j)=\{i,j\}$\\
$Fv(R_i v_{k_0}\cdots v_{k_{n_i-1}})=\{{k_0},\cdots {k_{n_i-1}}\}$ \> $Fv(\neg \phi)=Fv(\phi)$\\
$Fv(( \phi \wedge \theta))=Fv(\phi) \cup Fv(\theta$) \> $Fv(( \phi \vee \theta))=Fv(\phi) \cup Fv(\theta)$\\
$Fv(\exists v_k \phi)=Fv( \phi)-\{k\}$\\
\end{tabbing}
\end{Definition}

\begin{Definition} 
A formula $\phi$ is a sentence if $Fv(\phi)=\emptyset$.

\end{Definition}

\subsection{Semantics of the De Morgan valued-logic}
\begin{Definition} 
$P$ is an $\mathfrak{M}$-valued relation on a non-empty set $A$ if $P\in \mathcal{P}(A^n)^M$ and $\{P^q:q\in M\}$ is a partition of $A^n$.

$\mathfrak{A}=\langle A, P_i\rangle_{i\in I}$ is a $\mathfrak{M}$-structure if $A\neq \emptyset$ and $P_i$ is a $\mathfrak{M}$-valued relation on a $A$ for all $i\in I$.
\end{Definition}

Let $\phi$ be a formula and $\mathfrak{A}$ be a $\mathfrak{M}$-structure. $\phi^\mathfrak{A}\in  \mathcal{P}(|\mathfrak{A}|^n)^M$ is defined as follows:
\begin{Definition} All operations in the following are in the full $\mathfrak{M}-\textrm{CSA}_{\alpha}$ with base $|\mathfrak{A}|$. The definition is inductive.
\

$((R_i v_{k_0}\cdots v_{k_{n_i-1}})^\mathfrak{A})^q=\{s\in|\mathfrak{A}|^\omega:P_i^qs_{k_0}\cdots s_{ k_{n_{i-1}}}\}$

$(v_j\thickapprox v_k)^\mathfrak{A}=D_{jk}$

$((t_p)^\mathfrak{A}=\boldsymbol{U}_p$

$(\phi\vee\theta)^\mathfrak{A}=\phi^\mathfrak{A}\sqcup\theta^\mathfrak{A}$

$(\phi\wedge\theta)^\mathfrak{A}=\phi^\mathfrak{A}\sqcap\theta^\mathfrak{A}$

$(\neg\phi)^\mathfrak{A}=\circleddash\phi^\mathfrak{A}$

$(\exists v_k\phi)^\mathfrak{A}=\boldsymbol{C}_k\phi^\mathfrak{A}$

$(\gamma_p\phi)^\mathfrak{A}= \boldsymbol{\Delta}_p\phi^\mathfrak{A}$

\end{Definition}

Since $\{\phi^\mathfrak{A}:\phi \textrm{ is a formula}\}$ is closed under all operations of the full $\mathfrak{M}-\textrm{CA}_{\alpha}$ with base $|\mathfrak{A}|$, it is the universe of a subalgebra $\mathfrak{S}$ of the full $\mathfrak{M}-\textrm{CSA}_{\alpha}$ with base $|\mathfrak{A}|$.

We will now motivate the definition of $(\exists v_k\phi)^\mathfrak{A}$. Let $\mathfrak{A}$ be an $\mathfrak{M}$-structure and let $\mathfrak{A}^*=(\mathfrak{A},x)_{x\in |\mathfrak{A}|}$. Enlarge the language with individual constants $\underline{x}$ for $x\in |\mathfrak{A}|$. Define $\phi\binom{\kappa}{\underline{x}}$ as the result of replacing all free occurrences of $v_\kappa$ with $\underline{x}$. It is tedious, but not very difficult, to show that $(\phi\binom{\kappa}{\underline{x}}^\mathfrak{A})^p=x\circ_\kappa(\phi^\mathfrak{A})^{p}$. Recall that in an $\mathfrak{M}$-cylindric set algebra, $(L\sqcup K)^p=\bigcup_{q+r=p}L^q\cap K^r$. It is not difficult to show that this generalizes to

\begin{center}
$(L_1\sqcup L_2 \cdots \sqcup L_n)^p=\bigcup_{q_1+q_2+\cdots q_n=p}L^{q_1}\cap L^{q_2}\cap\cdots L^ {q_n}$
\end{center}
In the infinite case, this generalizes to
\begin{center}
 
$(\bigsqcup_{i\in I}L_i)^p=\bigcup_{j\in M^I,\hat{j}=p}\bigcap_{i\in I}L^{j_i}$

\end{center}
where $\hat{j}=sup\{j_i:i\in I\}$.

We regard $\exists v_\kappa \phi$ as the infinite disjunction of $\phi\binom{\kappa}{\underline{x}}$ for $x\in |\mathfrak{A}|$. Let $\bigvee_{x\in |\mathfrak{A}|}\phi\binom{\kappa}{\underline{x}}$ be this infinite disjunction. We then have

\begin{eqnarray*}
((\exists v_\kappa\phi)^\mathfrak{A})^p&=&\bigvee_{x\in |\mathfrak{A}|}(\phi\binom{\kappa}{\underline{x}}^\mathfrak{A})^p\\
&=&\bigcup_{\hat{j}=p}\bigcap_{x\in|\mathfrak{A}|}\phi\binom{\kappa}{\underline{x}}^\mathfrak{A})^{j_x} \textrm{	by above generalization}\\
&=&\bigcup_{\hat{j}=p}\bigcap_{x\in|\mathfrak{A}|}(x\circ_\kappa(\phi^\mathfrak{A}))^{j_x} \textrm{	}\\
&=&(\boldsymbol{C}_k\phi^\mathfrak{A})^p\\
\end{eqnarray*}

\begin{Theorem}
Let $\mathfrak{A}$ be a $\mathfrak{M}$-structure. If $k\notin Fv(\phi)$, then $(\exists v_k\phi)^\mathfrak{A}=\phi^\mathfrak{A}$.

\begin{proof}
It is straight forward to prove by induction on formulas that if $s,t\in|\mathfrak{A}|^\omega$ and $s\lceil Fv(\phi)=t\lceil Fv(\phi)$, then $s\in (\phi^\mathfrak{A})^p$ if and only if $t\in (\phi^\mathfrak{A})^p$ for all $p\in M$. Using the fact that $s\lceil Fv\phi=s\binom{\kappa}{x}\lceil Fv\phi$ if $\kappa \notin Fv\phi$, it is easy to show that $((\exists v_k\phi)^\mathfrak{A})^p=(\phi^\mathfrak{A})^p$ for all $p\in M$.

\end{proof}

\end{Theorem}

In the following, $\mathfrak{M}^*=\langle M,+,\cdot,-,0,1,\delta_p\rangle_{p\in M}$ where $\mathfrak{M}=\langle M,+,\cdot,-,0,1\rangle$ is a DeMorgan algebra and $\delta_p$ is defined by
\begin{displaymath}
\delta_pq=\left \{ \begin{array}{ll}
0 & \textrm{if $p\ne q $}\\
1 &  \textrm{if $p=q$}\end{array} \right .
\end{displaymath}
for $p,q \in M$.

\begin{Definition}

$\exists v_k\phi$, $v_j\thickapprox v_k$, and $R_i v_{k_0}\cdots v_{k_{n_i-1}}$ are prime formulas. Let $\mathbb{P}$ be the set of prime formulas and let $\phi$ be a formula and $v\in M^\mathbb{P}$. $T(\phi, v)\in M$ is defined as follows:\\
$T(\phi,v)=v(\phi) \textrm{ if } \phi \in \mathbb{P}$\\
$T(\phi\vee \theta,v)=T(\phi,v)+T(\theta,v)$\\
$T(\phi\wedge \theta,v)=T(\phi,v)\cdot T(\theta,v)$\\
$T(\neg\phi,v)=-T(\phi,v)$\\
$T(t_p,v)=p$\\
$T(\gamma_p\phi,v)=\delta_pT(\phi,v)$\\
$\phi$ is a tautology if $T(\phi,v)=1$ for all $v\in M^\mathbb{P}$.
\end{Definition}

\begin{Definition}

$\phi$ is true in a $\mathfrak{M}$-structure $\mathfrak{A}$ if $\phi^\mathfrak{A}=\boldsymbol{U}^1$ in the $\mathfrak{M}-\textrm{CSA}_{\alpha}$ $\mathfrak{S}$.\\
$\phi$ is a validity if $\phi$ is true in every $\mathfrak{M}$-structure $\mathfrak{A}$.\\
$\mathfrak{A}$ is a model of a set $\Sigma$ of formulas if $\phi$ is true in $\mathfrak{A}$ for all $\phi\in\Sigma$.\\
$\Sigma\models\phi$ if $\phi$ is true in every model of $\Sigma$.
\end{Definition}

\begin{Definition}
Let $Q\subseteq M$.
$\phi$ is $Q$-true in a $\mathfrak{M}$-structure $\mathfrak{A}$ if $\bigcup_{p\in Q}(\phi^\mathfrak{A})^p=\boldsymbol{U}^1$ in the $\mathfrak{M}-\textrm{CSA}_{\alpha}$ $\mathfrak{S}$.\\
$\phi$ is $Q$-validity if $\phi$ is $Q$-true in every $\mathfrak{M}$-structure $\mathfrak{A}$.\\
$\mathfrak{A}$ is a $Q$-model of $\Sigma$ if $\phi$ is $Q$-true in $\mathfrak{A}$ for all $\phi\in \Sigma$.\\
$\Sigma\models_Q\phi$ if $\phi$ is $Q$-true in every $Q$-model of $\Sigma$.
\end{Definition}

Since $\vee$ is commutative and associative in the sense that $(\phi\vee\theta)^\mathfrak{A}=\phi^\mathfrak{A}\sqcup\theta^\mathfrak{A}=\theta^\mathfrak{A}
 \sqcup\phi^\mathfrak{A}=(\theta\vee\phi)^\mathfrak{A}$ and $((\phi\vee\theta)\vee\psi))^\mathfrak{A}=(\phi^\mathfrak{A}\sqcup\theta^\mathfrak{A})\sqcup\psi^\mathfrak{A}=\phi^\mathfrak{A}\sqcup(\theta^\mathfrak{A}\sqcup\psi^\mathfrak{A})=(\phi\vee(\theta\vee\psi))^\mathfrak{A}$, we can generalize of $\vee$ to a finite disjunction $\bigvee_{i=I}\phi_i$ where $I$ is a finite set. In this case we obtain $(\bigvee_{i=I}\phi_i)^\mathfrak{A}=\bigsqcup_{i\in I}\phi_i^\mathfrak{A}$.

 \begin{Definition}
 Let $Q\subseteq M$. Then $\phi^Q=\bigvee_{p\in Q}\gamma_p\phi$. If $\Sigma$ is a set of formulas, then $\Sigma^Q=\{\phi^Q:\phi \in\Sigma\}$.
\end{Definition}

\begin{Theorem}
\label{Q true}
$\phi^Q$ is true in $\mathfrak{A}$ if and only if $\phi$ is $Q$-true in $\mathfrak{A}$.
\begin{proof}
The following computations take place in the $\mathfrak{M-}\textrm{CSA}_{\alpha}$ $ \mathfrak{S}$. Since $((\phi^Q)^\mathfrak{A})^q=\emptyset$ if $q\neq 0,1$ we can use Theorem \ref{insert}. First note that 
 $((\phi^Q)^\mathfrak{A})^1=((\bigvee_{q\in Q}\gamma_p\phi)^\mathfrak{A})^1=(\bigsqcup_{p\in Q}(\gamma_p\phi)^\mathfrak{A})^1=\bigcup_{p\in Q}((\gamma_p\phi)^\mathfrak{A})^1=\bigcup_{p\in Q}(\Delta_p\phi^\mathfrak{A})^1=\bigcup_{p\in Q}(\phi^\mathfrak{A})^p$. We then have
 
 \begin{eqnarray*}
 (\phi^Q)^\mathfrak{A}  \textrm{is true in } \mathfrak{A}  &\textrm{if and only if}&  ((\phi^Q)^\mathfrak{A})^1=\boldsymbol{U}^1\\
 &\textrm{if and only if}& \bigcup_{p\in Q}(\phi^\mathfrak{A})^p=\boldsymbol{U}^1\\
 &\textrm{if and only if}& \phi \textrm{ is } Q\textrm{-true in }\mathfrak{A}
  \end{eqnarray*}

\end{proof}
\end{Theorem}

\begin{Theorem}
\label{Q model}
$\mathfrak{A}$ is a model of $\Sigma^Q$ if and only if $\mathfrak{A}$ is a $Q$-model of $\Sigma$.
\begin{proof}
The proof follows directly from Theorem \ref{Q true}
\end{proof}

\end{Theorem}

\begin{Theorem}
\label{Q models}
$\Sigma\models_Q\phi$ if and only if $\Sigma^Q\models\phi^Q$.
\begin{proof}
The proof follows directly from Theorem \ref{Q model}
\end{proof}
\end{Theorem}

\begin{Theorem}
If $\phi$ is a tautology, then $\phi$ is a validity.
\begin{proof}
Let $\mathfrak{A}$ be a $\mathfrak{M}$-structure. For $s\in |\mathfrak{A}|^\omega$ and $\phi$ a formula, let $v_s\in U^\mathbb{P}$ be defined by $v_s(\phi)=p$ where $s\in (\phi^\mathfrak{A})^p$. Note that $v_s(\phi)=p$ is well defined since $\phi^\mathfrak{A}$ partitions $|\mathfrak{A}|^\omega$. Therefore, $v_s(\phi)=p$ if and only if $s\in (\phi^\mathfrak{A})^p$. By induction on formulas it is easily verified that:
\begin{center}
$T(\phi,v_s)=p$ if and only if $s\in (\phi^\mathfrak{A})^p$.\\
\end{center}
Let $\phi$ be a tautology. Therefore, $T(\phi,v)=1$ for all $v\in M^\mathbb{P}$. In particular, $T(\phi,v_s)=1$ for all $s\in |\mathfrak{A}|^\omega$ and hence, $s\in ((\phi)^\mathfrak{a})^1$ for all $s\in |\mathfrak{A}|^\omega$. It follows that $(\phi^\mathfrak{a})^1=|\mathfrak{A}|^\omega$ and $\phi^\mathfrak{A}=\boldsymbol{U}^1$ in the $\mathfrak{M}-\textrm{CSA}_{\alpha}$ $\mathfrak{S}$ and we have $\phi$ is true in $\mathfrak{A}$ and hence $\phi$ is a validity.
\end{proof}
\end{Theorem}

\begin{Definition}
In a $\mathfrak{M}-\textrm{CA}_{\alpha}$, $a\ |\!\!\!\Rightarrow b=\prod_{r\in M}(\delta_r b \boldsymbol{\mapsto} \sum_{q\le r}\delta_q a)$
\end{Definition}

\begin{Definition}
Let $\theta$ and $\phi$ be formulas. $\theta\rightarrow\phi=\neg\theta \vee \phi$, $\theta \Rightarrow \phi=\bigwedge_{r\in M}(\gamma_r\phi\to \bigvee_{q\le r}\gamma_q\theta)$, $\forall v_\kappa\phi=\neg\exists v_\kappa\neg\phi$, and $\theta\Leftrightarrow\phi=(\theta\Rightarrow\phi)\wedge(\phi\Rightarrow\theta)$ where $\bigwedge$ is generalized $\wedge$ and $\bigvee$ is generalized $\vee$.
\end{Definition}

Note that $(\theta\rightarrow\phi)^\mathfrak{A}=\theta^\mathfrak{A}\mapsto\phi^\mathfrak{A}$

We then have in $\mathfrak{S}$
\begin{displaymath}
(\theta\Rightarrow\phi)^\mathfrak{A}=(\bigwedge_{r\in M}(\gamma_r\phi\to \bigvee_{q\le r}\gamma_q\theta)^\mathfrak{A}=\prod_{r\in M}(\delta_r\phi^\mathfrak{A}\mapsto \sum_{q\le r}\delta_q\theta^\mathfrak{A})=\theta^\mathfrak{A}\ |\!\!\!\Rightarrow\phi^\mathfrak{A}
\end{displaymath}
Note that  $((\theta\Rightarrow \phi)^\mathfrak{A})^p=0$ if $p\ne 0,1$.

\begin{Theorem}
If $v\in M^\mathbb{P}$, then
\begin{displaymath}
T(\theta\Rightarrow \phi,v)=\left \{ \begin{array}{ll}
1 & \textrm{if $T(\theta,v)\le T(\phi,v)$}\\
0& \textrm{if $T(\theta,v)\nleq T(\phi,v)$}\end{array} \right .
\end{displaymath}
\begin{proof}

The following computations take place in $\mathfrak{M}^*$.
\begin{eqnarray*}
T(\theta\Rightarrow \phi,v)&=&T(\bigwedge_{r\in M}(\gamma_r\phi\rightarrow\bigvee_{q\le r}\gamma_q\theta),v)\\
&=&T(\bigwedge_{r\in M}(\neg\gamma_r\phi\vee\bigvee_{q\le r}\gamma_q\phi,v)\\
&=&\prod_{r\in M}T(\neg\gamma_r\phi\vee\bigvee_{q\le r}\gamma_q\theta,v)\\
&=&\prod_{r\in M}(T(\neg\gamma_r\phi,v)+T(\bigvee_{q\le r}\gamma_q\theta,v))\\
&=&\prod_{r\in M}(-\delta_rT(\phi,v)+\sum_{q\le r}\delta_qT(\theta,v))\\
&=&-\delta_{T(\phi,v)}T(\phi,v)+\sum_{q\le T(\phi,v)}\delta_qT(\theta,v)\\
&=&-1+\sum_{q\le T(\phi,v)}\delta_qT(\theta,v)\\
&=&\sum_{q\le T(\phi,v)}\delta_qT(\theta,v)
\end{eqnarray*}
\begin{displaymath}
=\left \{ \begin{array}{ll}
1 & \textrm{if $T(\theta,v)\le T(\phi,v)$}\\
0& \textrm{if $T(\theta,v)\nleq T(\phi,v)$}\end{array} \right .
\end{displaymath}

\end{proof}
\end{Theorem}

\begin{Theorem}
If $v\in M^\mathbb{P}$, then
\begin{displaymath}
T(\theta\Leftrightarrow \phi,v)=\left \{ \begin{array}{ll}
1 & \textrm{if $T(\theta,v)= T(\phi,v)$}\\
0& \textrm{if $T(\theta,v)\ne T(\phi,v)$}\end{array} \right .
\end{displaymath}
\begin{proof}
Note that $T(\theta\Leftrightarrow\phi,v)=T(\theta\Rightarrow\phi,v)\cdot T(\phi\Rightarrow\theta,v)$. If $v(\theta)=v(\phi)$, then $v(\theta)\le v(\phi)$ and $v(\phi)\le v(\theta)$. Therefore, $T(\theta\Rightarrow\phi,v)=1$ and $T(\phi\Rightarrow\theta,v)=1$ and hence, $T(\phi\Leftrightarrow\theta,v)=1$.\
If $v(\theta)\ne v(\phi)$, then $v(\theta)\nleq v(\phi)$ or $v(\phi)\nleq v(\theta)$. Therefore, $T(\theta\Rightarrow\phi,v)=0$ or $T(\phi\Rightarrow\theta,v)=0$ and hence, $T(\phi\Leftrightarrow\theta,v)=0$.

\end{proof}
\end{Theorem}

From the above, it follows that in $\mathfrak{M}^*$, $T(\theta\Rightarrow \phi,v)\in \{0,1\}$ and $T(\theta\Leftrightarrow\phi,v)\in \{0,1\}$.

\begin{Definition}
If $\phi$ is a formula, then $\Gamma\phi=\gamma_0\phi\vee\gamma_1\phi$.
\end{Definition}

Note that
\begin{displaymath}
T(\Gamma\phi,v)=\left \{ \begin{array}{ll}
1 & \textrm{if $T(\phi,v)\in\{0,1\}$}\\
0& \textrm{if $T(\phi,v)\notin\{0,1\}$}\end{array} \right .
\end{displaymath}

\begin{Theorem}
\label{Tautologies}
The following are tautologies where $Q\subseteq M$:\\
\begin{tabbing}
xxxxxxxxxxxxxxxxxxxxxxxxxxxxxxxxxxxx\= \kill
1. $(\phi\Rightarrow\theta)\rightarrow((\theta \Rightarrow\psi) \rightarrow(\phi \Rightarrow\psi))$\>2. $(\phi\Leftrightarrow\theta) \rightarrow(\phi \Rightarrow\theta)$\\
3. $(\phi\Leftrightarrow\theta) \rightarrow(\theta\Rightarrow \phi)$\>4. $(\theta \Rightarrow\phi) \rightarrow((\phi \Rightarrow\theta) \rightarrow(\phi \Leftrightarrow\theta)$\\
5. $\Gamma\phi\rightarrow(\Gamma(\phi\rightarrow(\theta\rightarrow\psi))\rightarrow(\Gamma(\phi\rightarrow\theta)\rightarrow((\phi\rightarrow(\theta\rightarrow\psi))\rightarrow((\phi\rightarrow\theta)\rightarrow(\phi\rightarrow\psi)))))$\\
6. $\Gamma\phi\rightarrow(\Gamma(\chi \Rightarrow(\phi \rightarrow\theta)) \rightarrow((\chi \Rightarrow(\phi \rightarrow\theta)) \rightarrow(\phi \rightarrow(\chi \Rightarrow\theta))))$\\
7. $\Gamma\phi\rightarrow[\Gamma(\phi \rightarrow(\chi \Rightarrow\theta)) \rightarrow[(\phi \rightarrow(\chi \Rightarrow\theta)) \rightarrow (\chi \Rightarrow (\phi \rightarrow \theta))]]$\\
8. $t_0^Q \rightarrow t_0$\>9.  $t_0 \rightarrow t_0^Q$\\
10. $\Gamma\phi \rightarrow\Gamma\neg\phi$\>11. $\Gamma\phi \rightarrow((\neg\phi \rightarrow t_0) \rightarrow \phi)$\\
12. $\phi \Rightarrow\phi$\>13  $\phi \Leftrightarrow\phi$\\
14. $(\phi \Leftrightarrow\theta) \rightarrow(\theta\Leftrightarrow\phi)$. \>15. $(\phi \Leftrightarrow\theta) \rightarrow((\theta\Leftrightarrow\psi) \rightarrow (\theta \Leftrightarrow \psi))$\\
16. $(\theta_1 \Leftrightarrow\psi_1) \rightarrow((\theta_2 \Leftrightarrow\psi_2) \rightarrow((\theta_1\vee\theta_2) \Leftrightarrow(\psi_1 \vee\psi_2))$\\
17. $(\theta_1 \Leftrightarrow\psi_1) \rightarrow((\theta_2 \Leftrightarrow\psi_2) \rightarrow((\theta_1\wedge\theta_2) \Leftrightarrow(\psi_1 \wedge\psi_2))$\\
18. $(\theta \Leftrightarrow\psi) \rightarrow(\neg\theta \Leftrightarrow\neg\psi)$\>19. $(\phi \Leftrightarrow\phi) \rightarrow(\gamma_p \Leftrightarrow \gamma_p\phi)$\\
20. $ (t_1 \Leftrightarrow t_0) \rightarrow t_0$\>21. $\Gamma\phi \rightarrow(\phi \rightarrow(\phi \Leftrightarrow t_1))$\\
22. $\gamma_1(\gamma_1\phi\vee\gamma_1\theta) \Leftrightarrow(\gamma_1\phi\vee\gamma_1\theta)$\>23. $\gamma_1(\gamma_1\phi\wedge\gamma_1\theta) \Leftrightarrow(\gamma_1\phi\wedge\gamma_1\theta)$\\

24. $\gamma_1(\neg\gamma_1\phi) \Leftrightarrow \neg\gamma_1\phi$\>25. $t_0 \Leftrightarrow\gamma_1t_0$\\
26. $t_1 \Leftrightarrow\gamma_1t_1$\>27.$ \gamma_p(\phi\vee\theta) \Leftrightarrow\bigvee_{q+r=p}(\gamma_q\phi\wedge\gamma_r\theta)$\\
28 $ \gamma_p(\phi\wedge\theta) \Leftrightarrow\bigvee_{q\cdot r=p}(\gamma_q\phi\wedge\gamma_r\theta)$\>29. $\gamma_p\phi\wedge\gamma_q\phi\Leftrightarrow t_0$ if $p\neq q$\\
30. $\bigvee_{p\in M}\gamma_p\phi \Leftrightarrow t_1$\>31. $\gamma_p(\neg\phi) \Leftrightarrow\gamma_{-p}\phi$\\
32. $\gamma_p t_p \Leftrightarrow t_1$\>33. $\gamma_p t_q \Leftrightarrow t_0$ if $p\neq q$\\
34. $\phi\vee\theta \Leftrightarrow\theta\vee \phi$\>35.  $\phi\wedge\theta \Leftrightarrow\theta\wedge \phi$\\
36. $\phi\vee(\theta \vee\psi) \Leftrightarrow(\phi\vee\theta) \vee\psi$\>37. $\phi \wedge(\theta \wedge\psi) \Leftrightarrow(\phi \wedge\theta) \wedge\psi$\\
38. $\phi\wedge(\theta\vee\psi) \Leftrightarrow(\phi\wedge\theta)\vee(\phi\wedge\psi)$\>39. $\phi \vee(\theta \wedge\psi) \Leftrightarrow(\phi \vee\theta) \wedge(\phi \vee\psi)$\\
40. $\phi\vee t_0 \Leftrightarrow\phi$\>41. $\phi\wedge t_1 \Leftrightarrow\phi$\\
42. $\gamma_1\phi\vee\neg\gamma_1\phi \Leftrightarrow t_1$\>43. $\gamma_1\phi \wedge\neg\gamma_1\phi \Leftrightarrow t_0$\\
44. $\gamma_1\gamma_q\phi \Leftrightarrow\gamma_q\phi$\>45. $\gamma_p\gamma_q\phi \Leftrightarrow t_0$ if $p\notin\{0,1\}$\\
46. $\gamma_0\gamma_q\phi \Leftrightarrow\neg\gamma_q\phi$\>47.  $\bigwedge_{p\in M}(\gamma_p\phi \Leftrightarrow\gamma_p\theta)\rightarrow(\phi \Leftrightarrow\theta)$\\
48. $\phi \Rightarrow(\phi\vee\theta)$\>49. $\Gamma\Gamma\phi$\\
50. $\Gamma\theta\rightarrow(\Gamma\phi \rightarrow(\theta \rightarrow(\phi \rightarrow(\theta\wedge\phi))))$\>51. $\Gamma\phi\rightarrow[\Gamma(\phi\rightarrow\psi)\rightarrow((\phi \rightarrow \psi) \rightarrow (\phi \rightarrow \Gamma\psi))]$\\

\end{tabbing}
\begin{proof}
As an example we prove 6. Note that $T(\psi \rightarrow \theta,v)=1$ if $T(\psi,v)=0$. Let $\psi=\Gamma\phi\rightarrow(\Gamma(\chi \Rightarrow(\phi \rightarrow\theta)) \rightarrow((\chi \Rightarrow(\phi \rightarrow\theta)) \rightarrow(\phi \rightarrow(\chi \Rightarrow\theta))))$. If $T(\phi,v)\notin \{0,1\}$, then $T(\Gamma\phi,v)=0$ and $T(\psi,v)=1$. So assume that $T(\phi,v)\in \{0,1\}$. Then $T(\psi,v)=\Gamma(\chi \Rightarrow(\phi \rightarrow\theta)) \rightarrow((\chi \Rightarrow(\phi \rightarrow\theta)) \rightarrow(\phi \rightarrow(\chi \Rightarrow\theta)))$. $T(\chi \Rightarrow(\phi \rightarrow\theta,v))\in\{0,1\}$ and hence, $T(\psi,v)=-T(\chi \Rightarrow(\phi \rightarrow\theta),v)+(-T(\phi,v) +T(\chi \Rightarrow \theta,v))$\\
If $T(\phi,v)=0$ then $T(\psi,v)=1$. So assume that $T(\phi,v)=1$. Then $T((\phi\rightarrow\theta,v)=T(\theta.v)$ and $T(\psi,v)=-T(\chi \Rightarrow(\phi \rightarrow\theta),v)+(T(\chi \Rightarrow \theta,v))$.\\
If $T(\chi,v)\leq T(\theta,v)$, then $T(\chi \Rightarrow\theta,v)=1$ and $T(\psi,v)=1$. So assume that $T(\chi,v)\nleq T(\theta,v)$. Therefore, $T(\chi \Rightarrow\theta,v)=0$ and we have $T(\psi,v)=-T(\chi \Rightarrow(\phi \rightarrow\theta,v)$. Since $T(\phi.v)=1$, we have $T(\phi \rightarrow\theta,v)=T(\theta,v)$ and $T(\chi,v)\nleq T(\phi \rightarrow\theta,v)$. Therefore, $T(\chi \Rightarrow(\phi \rightarrow\theta,v))=0$ and $T\psi,v)=1$.

\end{proof}
\end{Theorem}

\begin{Definition}
In a $\mathfrak{M}-\textrm{CA}_{\alpha}$, $s^\kappa_\lambda a=c_\kappa(d_{\kappa\lambda}\cdot a)$ and $q_\kappa a=-c_\kappa -a$. In a $\mathfrak{M}-\textrm{CSA}_{\alpha}$, $\boldsymbol{S}^\kappa_\lambda L=\boldsymbol{C}_\kappa(\boldsymbol{D}_{\kappa\lambda}\sqcap L)$ and $\boldsymbol{Q}_\kappa L= \boldsymbol{-}\boldsymbol{C}_\kappa\boldsymbol{-}L$. In our $\mathfrak{M}$-Valued Logic, $S^k_l\phi=\exists v_k(v_k\thickapprox v_l\wedge \phi)$.
\end{Definition}

It is easily shown that in the $\mathfrak{M}-\textrm{CSA}_{\alpha}$, $\mathfrak{S}$, $(S^k_l\phi)^\mathfrak{A}=\boldsymbol{S}^k_l\phi^\mathfrak{A}$ and $(\forall v_\kappa\phi)^\mathfrak{A}=\boldsymbol{Q}_\kappa\phi^\mathfrak{A}$.

In the following theorem we will need the following easy to show fact:

In an $\mathfrak{M}-\textrm{CSA}_{\alpha}$, if $L^q=0$ for $q\notin \{0,1\}$, then 
\begin{displaymath}
(L\sqcap K)^p=\left \{ \begin{array}{ll}
L^0\cup K^0 & \textrm{if $p=0$}\\
L^1\cap K^p &  \textrm{if $p\neq 0$}\end{array} \right .
\end{displaymath}

\begin{Theorem}
\label{substitution}
In a $\mathfrak{M}-\textrm{CSA}_{\alpha}$ with base $V$\\
	(a) $(\boldsymbol{S}^\kappa_\lambda L)^p=\bigcup_{j\in M^V, \hat{j}=p}(\bigcap_{x\in V,j_x=0}x\circ_\kappa( (\boldsymbol{D}_{\kappa\lambda})^0\cup L^0)\cap \bigcap_{x\in V,j_x\neq0}x\circ_\kappa((\boldsymbol{D}_{\kappa\lambda})^1\cap L^{j_x}))$\\
	(b) If $\kappa\neq\lambda$ and $t\in V^\alpha$, then  $t\in(\boldsymbol{S}^\kappa_\lambda L)^p$ if and only if $t\binom{\kappa}{t_\lambda}\in L^p$.\\
	(c) If $t\in V^\alpha$, then  $t\binom{\kappa}{t_\lambda}\in ((R_iv_{j_0}\cdots v_{j_{n_i-1}})^\mathfrak{A})^p$ if and only if\\
\begin{displaymath}
\left \{ \begin{array}{ll}
t\in ((R_iv_{j_0}\cdots v_{j_{l-1}}v_\lambda v_{j_{l+1}}\cdots v_{j_{n_i-1}})^\mathfrak{A})^p & \textrm{if $k=j_l$}\\
t\in ((R_iv_{j_0}\cdots  v_{j_{n_i-1}})^\mathfrak{A})^p  &  \textrm{if $k\neq j_l, l\neq 0,\cdots n_i-1$}\end{array} \right .
\end{displaymath}\\
(d) $(S^{k_0}_{j_0}\cdots S^{k_{n_i-1}}_{j_{n_i-1}}S^0_{k_0}\cdots S^{n_i-1}_{k_{n_i-1}}R_iv_0\cdots v_{n_i-1})^\mathfrak{A}=(R_iv_{j_0}\cdots v_{j_{n_i-1}})^\mathfrak{A}$ where $k_m\ne k_l$ for $m\ne l$ and $k_m\notin \{0,\cdots n-1,j_0\cdots j_{n_i-1}\}$ for $m<n_i$.

\begin{proof}
(a) This is easily shown using the above fact.\\
(b) Let $t\in(\boldsymbol{S}^\kappa_\lambda L)^p$. By (a), there is a $j\in M^V$ with, $\hat{j}=p$ such that if $j_x=0$, then $t\in x\circ_\kappa (\boldsymbol{D}_{\kappa\lambda})^0\cup L^0$, and if $j_x\neq 0$, then $t\in x\circ_\kappa(\boldsymbol{D}_{\kappa\lambda})^1\cap L^{j_x}$. Therefore, if $j_x=0$, then $t\binom{\kappa}{x}\in(\boldsymbol{D}_{\kappa\lambda})^0$ or $t\binom{\kappa}{x}\in L^0$ and if $j_x\neq 0$, then $t\binom{\kappa}{x}\in(\boldsymbol{D}_{\kappa\lambda})^1$ and $t\binom{\kappa}{x}\in L^{j_x}$. Therefore, we have\\
1. If $j_x=0$, then $x\neq t_\lambda$ or $t\binom{\kappa}{x}\in L^0$ and \\
2. If $j_x\neq 0$, then $x=t_\lambda$ and $t\binom{\kappa}{x}\in L^{j_x}$.\\
By 2, for all $x\neq t_\lambda$, $j_x=0$ and hence, $\hat{j}=j_{t_\lambda}$ and $j_{t_\lambda}=p$.\\
If $p=0$, then $j_{t_\lambda}=0$ and by 1 we have $t_\lambda\neq t_ \lambda$ or $t\binom{\kappa}{t_\lambda}\in L^0$ and hence, $t\binom{\kappa}{t_\lambda}\in L^0=L^p$.\\
If $p\neq 0$, then $j_{t_\lambda}\neq 0$ and by 2 we have $t_\lambda=t_\lambda$ and $t\binom{\kappa}{t_\lambda}\in L^{j_{t_\lambda}}$ and hence, $t\binom{\kappa}{t_\lambda}\in L^{j_{t_\lambda}}=L^p$.\\
Let $\binom{\kappa}{t_\lambda}\in L^p$. Define $j\in M^V$ by
\begin{displaymath}
j_x=\left \{ \begin{array}{ll}
0 & \textrm{if $x\neq t_\lambda$}\\
p  &  \textrm{if $x=t_\lambda$}\end{array} \right .
\end{displaymath}
Then $\hat{j}=p$. Note that if $p=0$, then $j_x=0$ for all $x\in V$.\\
Case 1. $p=0$. Then by (a), $(\boldsymbol{S}^\kappa_\lambda L)^p=\bigcup_{j\in M^V, \hat{j}=p}\bigcap_{x\in V,j_x=0}x\circ_\kappa( (\boldsymbol{D}_{\kappa\lambda})^0\cup L^0)$. If $x= t_\lambda$, then $t\binom{\kappa}{x}=t\binom{\kappa}{t_\lambda}\in L^p=L^0$. If $x\neq t_\lambda$, then $t\binom{\kappa}{x}\in (\boldsymbol{D}_{\kappa\lambda})^0$. In either case, $t\binom{\kappa}{x}\in (\boldsymbol{D}_{\kappa\lambda})^0\cup L^0$ and hence, $t\in x\circ_\kappa ( (\boldsymbol{D}_{\kappa\lambda})^0\cup L^0)$. Therefore, $t\in(\boldsymbol{S}^\kappa_\lambda L)^0$.\\
Case 2. $p\neq 0$. If $j_x=0$, then $x\neq t_\lambda$ and hence, $t\binom{\kappa}{x}\in (\boldsymbol{D}_{\kappa\lambda})^0$ and $t\binom{\kappa}{x}\in (\boldsymbol{D}_{\kappa\lambda})^0\cup L^0$. If $j_x\neq 0$, then $x=t_\lambda$ and $j_x=p$. Therefore, $t\binom{\kappa}{x}\in (\boldsymbol{D}_{\kappa\lambda})^1$ and $t\binom{\kappa}{x}=t\binom{\kappa}{t_\lambda}\in L^p=L^{j_x}$. Therefore, $t\binom{\kappa}{x}\in (\boldsymbol{D}_{\kappa\lambda})^1\cap L^{j_x}$ and $t\in x\circ_\kappa((\boldsymbol{D}_{\kappa\lambda})^1\cap L^{j_x}))$. Therefore, $t\in(\boldsymbol{S}^\kappa_\lambda L)^p$.\\
(c) Follows easily from the definitions.\\
(d) Follows easily using repeated applications of (b) followed by repeated applications of (c).
\end{proof}
\end{Theorem}

\begin{Theorem}
\label{less than true}
$\theta\Rightarrow\phi$ is a true in $\mathfrak{A}$ if and only if $\theta^\mathfrak{A}\leq\phi^\mathfrak{A}$.
\begin{proof}
It is easy to prove this using Theorems \ref{insert} and \ref{less than}.
\end{proof}
\end{Theorem}

\begin{Theorem}
\label{equal true}
$\theta\Leftrightarrow\phi$ is true in $\mathfrak{A}$ if and only if $\theta^\mathfrak{A}=\phi^\mathfrak{A}$.
\begin{proof}
The following computation takes place in $\mathfrak{S}$.\\
$((\theta\Leftrightarrow\phi)^\mathfrak{A})^1=((\theta\Rightarrow\phi\wedge\phi\Rightarrow\theta)^\mathfrak{A})^1 =((\theta\Rightarrow\phi)^\mathfrak{A}\sqcap (\phi\Rightarrow\theta)^\mathfrak{A})^1= \bigcup_{q\cdot r=1}[((\phi\Rightarrow\theta)^\mathfrak{A})^q\cap((\theta\Rightarrow\phi)^\mathfrak{A})^r]= \bigcup_{q\cdot r=1}[(\phi^\mathfrak{A}|\!\!\!\Rightarrow \theta^\mathfrak{A})^q\cap(\theta^\mathfrak{A}|\!\!\!\Rightarrow \phi^\mathfrak{A})^r]= (\phi^\mathfrak{A}|\!\!\!\Rightarrow \theta^\mathfrak{A})^1\cap (\theta^\mathfrak{A}|\!\!\!\Rightarrow \phi^\mathfrak{A})^1$\\
Therefore, we have
\begin{eqnarray*}
\theta\Leftrightarrow\phi \textrm{ is true in } \mathfrak{A}&\textrm{if and only if }& ((\theta\Leftrightarrow\phi)^\mathfrak{A})^1 =|\mathfrak{A}|^\omega\\
&\textrm{if and only if }&  (\phi^\mathfrak{A}|\!\!\!\Rightarrow \theta^\mathfrak{A})^1\cap (\theta^\mathfrak{A}|\!\!\!\Rightarrow \phi^\mathfrak{A})^1=|\mathfrak{A}|^\omega\\
&\textrm{if and only if }&(\phi^\mathfrak{A}|\!\!\!\Rightarrow \theta^\mathfrak{A})^1=|\mathfrak{A}|^\omega \textrm{and } (\theta^\mathfrak{A}|\!\!\!\Rightarrow \phi^\mathfrak{A})^1=|\mathfrak{A}|^\omega\\
&\textrm{if and only if }&((\phi\Rightarrow\theta)^\mathfrak{A})^1=|\mathfrak{A}|^\omega \textrm{ and } ((\theta\Rightarrow\phi)^\mathfrak{A})^1=|\mathfrak{A}|^\omega\\
&\textrm{if and only if }& \phi\Rightarrow\theta \textrm{ is true in }  \mathfrak{A} \textrm{ and } \theta\Rightarrow\phi \textrm{ is true in }  \mathfrak{A}\\
&\textrm{if and only if }& \phi^\mathfrak{A}\le \theta^\mathfrak{A} \textrm{ and }\theta^\mathfrak{A}\le \phi^\mathfrak{A} \\
&\textrm{if and only if }& \phi^\mathfrak{A}= \theta^\mathfrak{A}
\end{eqnarray*} 

\end{proof}
\end{Theorem}

In a $\mathfrak{M}-\textrm{CSA}_{\alpha}$ with base $U$, $\boldsymbol{Q}_\kappa A=\bigcup_{w\in U^M,\check{w}=p}\bigcap_{x\in U}x\circ_\kappa A^{w_x}$ where $\check{w}=\prod_{x\in U}w_x$.

\begin{Theorem}
\label{validities}
The following are validities:\\
\begin{tabbing}
xxxxxxxxxxxxxxxxxxxxxxxxxxxxxxxxxxxx\= \kill
1. $\Gamma \forall v_k\phi\rightarrow(\forall u_k\phi\rightarrow\phi)$\>2. $\gamma_1\exists v_k\gamma_1\phi\Leftrightarrow\exists v_k\gamma_1\phi$\\
3. $\gamma_1(v_k\thickapprox v_l)\Leftrightarrow v_k\thickapprox v_l$\>4. $\gamma_p(v_k \thickapprox v_l) \Leftrightarrow t_0$ if $p\notin \{0,1\}$\\
5. $\gamma_0(v_k \thickapprox v_l)\Leftrightarrow\neg v_k \thickapprox v_l $\>6.  $\exists v_kt_0 \Leftrightarrow t_0$\\
7. $\phi\vee\exists v_k\phi \Leftrightarrow\exists v_k\phi$\>8. $\exists v_k(\phi\wedge\exists v_k\theta) \Leftrightarrow\exists v_k\phi\wedge\exists v_k\theta$\\
9. $\exists v_k \exists v_l \Leftrightarrow \exists v_l \exists v_k\phi$\>10. $v_k \thickapprox v_k \Leftrightarrow t_1$\\
11. $v_l \thickapprox v_m \Leftrightarrow\exists v_k(v_l \thickapprox v_k\wedge v_k \thickapprox v_m)$ if $k\notin\{l,m\}$\\
12. $\exists v_k(v_k \thickapprox v_l\wedge\gamma_1\phi)\wedge\exists v_\kappa(v_k \thickapprox v_l\wedge\neg\gamma_1\phi) \Leftrightarrow t_0$ if $k\ne l$\\
13. $\gamma_p(\exists v_k\phi) \Leftrightarrow\bigvee_{A\subseteq M,supA=p}\bigwedge_{q\in A}\exists v_k\delta_{w_q}\phi\wedge\neg\bigvee_{A\subseteq M,supA>p}\bigwedge_{q\in A} \exists v_k\gamma_p{w_q}\phi$\\
14. $R_iv_{j_0}\cdots v_{j_{{n_i}-1}} \Leftrightarrow S^{k_0}_{j_0}\cdots S^{k_{n_i-1}}_{j_{n_i-1}}S^0_{k_0}\cdots S^{n_{i-1}}_{k_{n_i-1}}R_iv_0\cdots v_{n_i-1}$ where $k_m\ne k_l$\\ for $m\ne l$ and $k_m\notin \{0,\cdots n-1,j_0\cdots j_{n_i-1}\}$ for $m<n_i$.\\
\end{tabbing}

\begin{proof}
Proof of 1: We use the following easily verified facts in a $\mathfrak{M}-\textrm{CSA}_{\alpha}$: $(\boldsymbol{\Delta}L)^1=L^1\cup L^0$ and $(\boldsymbol{\Delta}L)^0=-(L^1\cup L^0)$ and hence, $(\boldsymbol{\Delta}L)^0=-(\boldsymbol{\Delta}L)^1$.\\
Let $\chi=\Gamma \forall v_k\phi\rightarrow(\forall u_k\phi\rightarrow\phi)$. Then in $\mathfrak{S}$, $(\chi^\mathfrak{A})^1=(\circleddash\boldsymbol{\Delta}\boldsymbol{Q}_k\phi^\mathfrak{A}\sqcup(\circleddash \boldsymbol{Q}_k\phi^\mathfrak{A} \sqcup\phi^\mathfrak{A}))^1$. We show that in a $\mathfrak{M}-\textrm{CSA}_{\alpha}$ with base $V$, $(\circleddash\boldsymbol{\Delta}\boldsymbol{Q}_k\L\sqcup(\circleddash \boldsymbol{Q}_k L \sqcup L))^1=V^\alpha$.\\

\begin{eqnarray*}
((\circleddash\boldsymbol{\Delta}\boldsymbol{Q}_k L \sqcup(\circleddash \boldsymbol{Q}_k L \sqcup L))^1&=&\bigcup_{r+t=1}((\circleddash\boldsymbol{\Delta} \boldsymbol{Q}_k L)^r\cap (\circleddash \boldsymbol{Q}_k L\sqcup L)^t)\\
&=&\bigcup_{r+t=1}((\boldsymbol{\Delta} \boldsymbol{Q}_k L)^{-r}\cap (\circleddash \boldsymbol{Q}_k L\sqcup L)^t)\\
&=&[\bigcup_{t\in M}((\boldsymbol{\Delta} \boldsymbol{Q}_k L)^0\cap (\circleddash \boldsymbol{Q}_k L\sqcup L)^t)]\cup[(\boldsymbol{\Delta} \boldsymbol{Q}_k L)^1\cap (\circleddash \boldsymbol{Q}_k L\sqcup L)^1]\\ 
& &\textrm{(breaking the union into 2 parts - one with $r=1$ and the other}\\
& & \textrm{with $r=0$. All other terms are $0$ since $(\boldsymbol{\Delta} \boldsymbol{Q}_k L)^r=0$ if $r\notin \{0,1\}$})\\
&=&[(\boldsymbol{\Delta} \boldsymbol{Q}_k L)^0\cap\bigcup_{t\in M} (\circleddash \boldsymbol{Q}_k L\sqcup L)^t]\cup[(\boldsymbol{\Delta} \boldsymbol{Q}_k L)^1\cap (\circleddash \boldsymbol{Q}_k L\sqcup L)^1]\\
&=&[(\boldsymbol{\Delta} \boldsymbol{Q}_k L)^0\cap V^\alpha]\cup[(\boldsymbol{\Delta} \boldsymbol{Q}_k L)^1\cap (\circleddash \boldsymbol{Q}_k L\sqcup L)^1]\\
&=&(\boldsymbol{\Delta} \boldsymbol{Q}_k L)^0\cup[(\boldsymbol{\Delta} \boldsymbol{Q}_k L)^1\cap (\circleddash \boldsymbol{Q}_k L\sqcup L)^1]\\
&=&(\boldsymbol{\Delta} \boldsymbol{Q}_k L)^0\cup[(\boldsymbol{\Delta} \boldsymbol{Q}_k L)^1\cap \bigcup_{r+t=1}(\circleddash \boldsymbol{Q}_k L)^r\cap L^t)]\\
&=&[(\boldsymbol{\Delta} \boldsymbol{Q}_k L)^0\cup(\boldsymbol{\Delta} \boldsymbol{Q}_k L)^1]\cap[(\boldsymbol{\Delta} \boldsymbol{Q}_k L)^0\cup \bigcup_{r+t=1}( \boldsymbol{Q}_k L)^{-r}\cap L^t]\\
&=&V^\alpha\cap[\boldsymbol{\Delta} \boldsymbol{Q}_k L)^0\cup \bigcup_{r+t=1}( \boldsymbol{Q}_k L)^{-r}\cap L^t]\\
&=&(\boldsymbol{\Delta} \boldsymbol{Q}_k L)^0\cup \bigcup_{r+t=1}( \boldsymbol{Q}_k L)^{-r}\cap L^t)\\
&=&-(\boldsymbol{\Delta} \boldsymbol{Q}_k L)^1\cup \bigcup_{r+t=1}( \boldsymbol{Q}_k L)^{-r}\cap L^t)
\end{eqnarray*}

It suffices to show that $(\boldsymbol{\Delta} \boldsymbol{Q}_k L)^1\subseteq \bigcup_{r+t=1}( \boldsymbol{Q}_k L)^{-r}\cap L^t)$; that is, $(\boldsymbol{Q}_k L)^1\cup (\boldsymbol{Q}_k L)^0\subseteq\bigcup_{r+t=1}( \boldsymbol{Q}_k L)^{-r}\cap L^t)$.\\
Let $s\in ( \boldsymbol{Q}_k L)^1=\bigcup_{\check{w}=1,w\in M^V}\bigcap_{x\in V}x\circ_k L^{w_x}=\bigcap_{x\in V}x\circ_k L^1$ since if $ \check{w}=1$ then $w_x=1$ for all $x\in V$. Therefore, $s\in x\circ_k L^1$ for all $x\in V$; that is, $s\binom{k}{x}\in L^1$ for all $x\in V$. Therefore,  $s=s\binom{k}{s_k}\in L^1$. Let $r=0$, $t=1$ and $j_x=1$ for all $x\in V$. Then $r+t=1$ and $\check{j}=1=-r$. Therefore, $s\in \bigcup_{\check{j}=-r}\bigcap_{x\in V}x\circ_k L^{j_x}$ and therefore, $s\in\bigcup_{r+t=1}(\boldsymbol{Q}_k L)^{-r}\cap L^t$\\
Let $s\in ( \boldsymbol{Q}_k L)^0=\bigcup_{\check{w}=0,w\in M^V}\bigcap_{x\in V}x\circ_k L^{w_x}$. There is a $w\in V^M$ such that $\check{w}=0$ and for all $x\in V$, $s\in x\circ_k L^{w_x}$. Therefore, $s\in s_k \circ_k L^{w_{s_k}}$ and hence, $s=s\binom{k}{s_k}\in L^{w_{s_{k}}}$. Let $t=w_{s_{k}}$ and $r=1$. Then $s\in L^t$ and $r+t=1$. Since $\check{w}=-r$, $s\in \bigcup_{\check{w}=-r}\bigcap_{x\in V}x\circ_k L^{w_x}=(\boldsymbol{Q}_k L)^{-r}$ and $s\in \bigcup_{r+t=1}( \boldsymbol{Q}_k L)^{-r}\cap L^t)$.\\
Proofs of 2 - 13 follow from Theorem \ref{equal true} and the appropriate axiom for $\mathfrak{M}-\textrm{CA}_{\alpha}$'s.\\
Proof of 14 follows from Theorem \ref{substitution} and Theorem \ref{equal true}.

\end{proof}

\end{Theorem}

\subsection{Proof Structure of $\mathfrak{M}$-Valued Logic}

In the following, we give a set of axioms and rules of inference that are needed to prove the completeness theorem. It is not the most elegant set of axioms and rules of inference but it allows us to prove the completeness theorem without having to prove theorems within the system.

All tautologies are axioms and, in addition, the validities of Theorem \ref{validities} are axioms.

Rules of infrerence:\\
1. Modus Ponens: From $\phi$ and $\phi \rightarrow\theta$, infer $\theta$.\\
2. $\Gamma$-rule: From $\phi$ infer $\Gamma\phi$.\\ 
3. $\exists$-rule: From $\theta \Rightarrow\phi$ infer $\exists v_k\theta \rightarrow\phi$ if $k\notin Fv\phi$.

\begin{Definition}
If $\Sigma$ is a set of formulas and $\phi$ is a formula then $\Sigma\vdash\phi$ means there is a finite sequence of formulas $\theta_1,\theta_2,\cdots \theta_n$ such that $\theta_n=\phi$. In addition, for $i\in \{1,2,\cdots,n\}$, $\theta_i \in\Sigma$ or $\theta_i$ is an axiom or\\
1. (Modus Ponens) $\theta_j=\theta_k\rightarrow\theta_i$ where $j,k<i$ or\\
2. ($\Gamma$-rule) $\theta_i=\Gamma\theta_j$ where $j<i$ or\\
3. ($\exists$-rule) $\theta_i=\exists v_k\phi\Rightarrow\psi$ where $\theta_j=\phi\Rightarrow\psi$, $k\notin Fv\psi$, and $j<k$
\end{Definition}

\begin{Definition}
It $Q\subseteq M$ then $\Sigma\vdash_Q\phi$ if $\Sigma^Q\vdash\phi^Q$.
\end{Definition}

\begin{Theorem} (Soundness Theorem)\\
(a) If $\Sigma\vdash\phi$, then $\Sigma\models\phi$.\\
(b) If $\Sigma\vdash_Q\phi$, then $\Sigma\models_Q\phi$.
\begin{proof}
(a) It suffices to show that for all $\mathfrak{M}$-structures $\mathfrak{A}$, that in $\mathfrak{S}$, the following hold:\\
1. If $\phi$ is an axiom, then $\phi^\mathfrak{A}= \boldsymbol{U}_1$.\\
2. If $(\theta \rightarrow\phi)^\mathfrak{A}= \boldsymbol{U}_1$ and $\theta^\mathfrak{A}= \boldsymbol{U}_1$ then $\phi^\mathfrak{A}= \boldsymbol{U}_1$.\\
3. If $\phi ^\mathfrak{A}= \boldsymbol{U}_1$, then $(\Gamma\phi)^\mathfrak{A}=\boldsymbol{U}_1$.\\
4. If $(\phi\Rightarrow\theta)^\mathfrak{A}=\boldsymbol{U}_1$ and $k\notin Fv\phi$, then $(\exists\  v_k\theta \Rightarrow\phi)^\mathfrak{A}=\boldsymbol{U}_1$\\
All of the above are easily verified.\\
(b) This follows from (a) and Theorem \ref{Q models}.
\end{proof}

\end{Theorem}

\begin{Theorem}
\label{Theorems}
\text{ }\\
(a) $\vdash\phi\Rightarrow\exists v_k\phi$.\\
(b) If $\Sigma\vdash\forall v_k\phi$, then $\Sigma\vdash\phi$.\\
(c) If $\Sigma\vdash\phi\Rightarrow\theta$, then $\Sigma\vdash\exists v_k\phi \Rightarrow\exists v_k\theta$.\\
(d) If $\Sigma\vdash\phi\Leftrightarrow\theta$, then $\Sigma\vdash\exists v_k\phi \Leftrightarrow\exists v_k\theta$.\\
(e) If $\kappa\notin Fv\phi$, then $\vdash\exists v_\kappa\phi \Leftrightarrow\phi$.\\
(f)  If $\Sigma\vdash\theta_i$ for $i\in \{1,2,\cdots ,n\}$ then $\Sigma\vdash\bigwedge_{i\in \{1,2,\cdots ,n\}}\theta_i$.

\begin{proof}\text{ }\\
(a) Here is a proof: $\phi\Rightarrow(\phi\vee\exists v_k\phi)$ (Tautology), $(\phi\vee\exists v_k\phi)\Leftrightarrow\exists v_k\phi$ (Axiom), $((\phi\vee\exists v_k\phi)\Leftrightarrow\exists v_k\phi) \rightarrow((\phi\vee\exists v_k\phi)\Rightarrow\exists v_k\phi)$ (Tautology), $(\phi\vee\exists v_k\phi)\Rightarrow\exists v_k\phi$ (MP), $[(\phi\Rightarrow(\phi\vee\exists v_k\phi))\rightarrow((\phi\vee\exists v_k\phi)\Rightarrow\exists v_k\phi)]\rightarrow(\phi\Rightarrow\exists v_k\phi)$ (Tautology)
Now use Modus Ponens twice to obtain $\phi\Rightarrow\exists v_k\phi$.\\
(b) Assume $\Sigma\vdash\forall v_k\phi$. By $\Gamma$-rule, $\Sigma\vdash\Gamma\forall v_k\phi$. Use the axiom $\Gamma\forall v_k\phi\rightarrow(\forall v_k\phi\rightarrow\phi)$ and Modus Ponens twice to obtain $\Sigma\vdash\phi$.\\
(c) Suppose that $\Sigma\vdash\phi\Rightarrow\theta$. Use part (a) $\phi\Rightarrow\exists v_k\theta$, the tautology $(\phi\Rightarrow\theta)\rightarrow((\theta\Rightarrow\exists v_k\theta)\rightarrow(\phi\Rightarrow\exists v_k\theta))$ and Modus Ponens twice to obtain $\phi\Rightarrow\exists v_k\theta$. Then use the $\exists$-rule to obtain $\Sigma\vdash\exists v_k\phi \Leftrightarrow\exists v_k\theta$.\\
(d) Suppose that $\Sigma\vdash\phi\Leftrightarrow\theta$. Use the tautology $(\phi\Leftrightarrow\theta)\rightarrow(\phi\Rightarrow\theta)$, Modus Ponens, and (c) to obtain $\Sigma\vdash\exists v_k\phi \Rightarrow\exists v_k\theta$. Use the tautology $(\phi\Leftrightarrow\theta)\rightarrow(\theta\Rightarrow\phi)$, Modus Ponens, and (b) to obtain $\Sigma\vdash\exists v_k\theta \Rightarrow\exists v_k\phi$. Now use the tautology $(\exists v_k\phi \Rightarrow\exists v_k\theta)\rightarrow((\exists v_k\theta \Rightarrow\exists v_k\phi)\rightarrow(\exists v_k\phi \Leftrightarrow\exists v_k\theta))$ and Modus Ponens twice to obtain $\Sigma\vdash\exists v_k\phi \Leftrightarrow\exists v_k\theta$.\\
(e) Use the following steps: $\phi\Rightarrow\phi$ (tautology), $\exists v_\kappa\phi \Rightarrow\phi$, ($\exists$-rule), $\phi \Rightarrow\exists v_\kappa\phi$ (Part a), $(\exists v_\kappa\phi \Rightarrow\phi)\rightarrow((\phi \Rightarrow\exists v_\kappa\phi) \rightarrow(\exists v_\kappa\phi \Leftrightarrow\phi$)) (tautology) and Modess Ponens twice to obtain $\vdash\exists v_\kappa\phi \Leftrightarrow\phi$.\\
(f) We sketch a proof for $n=2$. The general proof is easily done by induction. Assume that  $\Sigma\vdash\theta$ and $\Sigma\vdash\phi$. Use the $\Gamma$-rule, the tautology $\Gamma\theta\rightarrow(\Gamma\phi \rightarrow(\theta \rightarrow(\phi \rightarrow(\theta\wedge\phi))))$ and Modus Ponens to obtain $\Sigma\vdash\theta\wedge\phi$.\\

\end{proof}
\end{Theorem}

\begin{Theorem} Deduction Theorem\\
Let $\Sigma$ be a set of formula, $\phi$ a sentence, and $\theta$ a formula. If $\Sigma\vdash\Gamma\phi$ and $\Sigma\cup\{\phi\}\vdash\theta$, then $\Sigma\vdash\phi\rightarrow\theta$.
\begin{proof}
Let $\chi_1,\chi_2,\cdots, \chi_m$ be a proof of $\Sigma\vdash\Gamma\phi$ and let $\psi_1, \psi_2,\cdots, \psi_n $ be a proof of $\Sigma\cup\{\phi\}\vdash\theta$. An outline of a proof of $\Sigma\vdash\phi\rightarrow\theta$ is $\chi_1,\chi_2,\cdots, \chi_m,\phi\rightarrow\psi_1, \phi\rightarrow\psi_2,\cdots, \phi\rightarrow\psi_n$. We will justify $\phi\rightarrow\psi_m$, $m\in\{1,\cdots,n\}$.\\
If $\psi_m$ is an axiom or $\psi_m\in\Sigma$, insert before $ \phi\rightarrow\psi_m$ the following: $\psi_m$, $\Gamma\psi_m$ ($\Gamma$-rule), $\Gamma\psi_m\rightarrow(\psi_m\rightarrow(\phi\rightarrow\psi_m))$ (Tautology), $\psi_m\rightarrow(\phi\rightarrow\psi_m)$, (MP), $\phi\rightarrow\psi_m$ (MP).\\
If $\psi_m=\phi$, insert the following: $\Gamma\phi\rightarrow(\phi\rightarrow\psi_m)$ (Tautology), $\phi\rightarrow\phi=\phi\rightarrow\psi_m$ (MP since $\Gamma\phi=\chi_l$)\\
$\exists$-rule: Assume that $\psi_m=\exists v_k\chi\Rightarrow\theta$ and $\psi_j=\chi\Rightarrow\theta$ where $j<m$ and $k\notin Fv\theta$ and $\phi\rightarrow\psi_j$ is already justified; that is $\phi\rightarrow(\chi\rightarrow\theta)$ is already justified. We will now justify $\phi\rightarrow \psi_m$.  Precede $\phi\rightarrow \psi_m$ with $\Gamma(\phi\rightarrow(\chi\Rightarrow\theta))$ ($\Gamma$-rule) and tautology $\Gamma\phi\rightarrow(\Gamma(\phi\rightarrow(\chi\Rightarrow\theta))\rightarrow((\phi\rightarrow(\chi\Rightarrow\theta))\rightarrow(\chi\Rightarrow(\phi\rightarrow\theta))))$. Use Modus Ponens three times to obtain $\chi\Rightarrow(\phi\rightarrow\theta)$. Since $k\notin Fv\theta$ and $\phi$ is a sentence, $k\notin Fv(\phi\rightarrow\theta)$. Insert $\exists v_k\chi\Rightarrow(\phi\rightarrow\theta)$ ($\exists$-rule) and then insert $\Gamma(\phi\rightarrow(\chi\Rightarrow\theta))$ ($\Gamma$-rule) and tautology $\Gamma\phi\rightarrow(\Gamma(\exists v_k \chi\Rightarrow(\phi\rightarrow\theta))\rightarrow((\exists v_k\chi\Rightarrow(\phi\rightarrow\theta))\rightarrow(\phi\rightarrow(\exists v_k\chi\Rightarrow\theta))))$. Now use Modus Ponens three times to obtain $\phi\rightarrow(\exists v_k\chi\Rightarrow\theta)$.\\
Modus Ponens: Assume that $\psi_m=\psi_j\rightarrow\psi_l$ where $m,j<l$. Assume that $\phi\rightarrow\psi_m=\phi\rightarrow(\psi_j\rightarrow\psi_l)$ and $\phi\rightarrow\psi_j$ have already been justified. We show that $\phi\rightarrow\psi_l$ is justified. Precede $\phi\rightarrow\psi_l$ with $\Gamma(\phi\rightarrow(\psi_j\rightarrow\psi_l))$ ($\Gamma$ -rule), $\Gamma(\phi\rightarrow\psi_j)$ ($\Gamma$ -rule), and the tautology $\Gamma\phi\rightarrow(\Gamma(\phi\rightarrow(\psi_j\rightarrow\psi_l))\rightarrow(\Gamma(\phi \rightarrow \psi_j) \rightarrow (\phi \rightarrow (\psi_j \rightarrow \psi_l) \rightarrow ((\phi \rightarrow \psi_j )\rightarrow (\phi \rightarrow \psi_l)))))$. Now use Modus Ponens five times to obtain $\phi \rightarrow\psi_l$.\\
$\Gamma$-rule: Assume that $\psi_j=\Gamma\psi_m$ where $m<j$ and $\phi \rightarrow\psi_m$ already justified. We show that $\phi \rightarrow\Gamma\psi_m$ is justified. Precede $\phi \rightarrow\Gamma\psi_m$ with $\Gamma(\phi \rightarrow\psi_m)$ ($\Gamma$-rule) and the tautology $\Gamma\phi \rightarrow(\Gamma(\phi \rightarrow\phi_m) \rightarrow((\phi \rightarrow\psi_m) \rightarrow(\phi \rightarrow\Gamma\Psi_m)))$. Now use Modus Ponens three times to obtain $\phi \rightarrow\Gamma\Psi_m$.
\end{proof}\
\end{Theorem}

\begin{Definition}
If $\Sigma$ is a set of formulas, then $\Sigma$ is consistent if it is not the case that $\Sigma\vdash t_0$\\
If $\Sigma$ is a set of formulas and $Q\subseteq M$, then $\Sigma$ is $Q$- consistent if it is not the case that $\Sigma\vdash_Q t_0$
\end{Definition}

\begin{Theorem}
\label{Q-consistent}
Let $\Sigma$ be a set of formulas, $Q\subseteq M$, and $0\notin Q$. $\Sigma$ is Q-consistent if and only if $\Sigma^Q$ is consistent.
\begin{proof}
This is easily seen using the tautologies $t_0^Q\rightarrow t_o$ and $t_0\rightarrow t_o^Q$.
\end{proof}
\end{Theorem}

\begin{Theorem}
\label{consistency}
Let $\phi$ be a sentence and $\Sigma$ a set of formulas. If $\Sigma\vdash \Gamma\phi$ and it is not the case that $\Sigma \vdash \phi$, then $\Sigma\cup\{\neg\phi\}$ is consistent.
\begin{proof}
Assume that $\Sigma\cup\{\neg\phi\}$ is inconsistent. Then $\Sigma\cup\{\neg\phi\}\vdash t_o$. Use the tautology $\Gamma\phi\rightarrow\Gamma\neg\phi$ and Modus Ponens to obtain $\Sigma\vdash \Gamma\neg\phi$. By the Deduction Theorem,  $\Sigma\vdash \neg\phi\rightarrow t_0$. Use the tautology $\Gamma\phi \rightarrow((\neg\phi \rightarrow t_o) \rightarrow\phi)$ and Modus Ponens twice to obtain $\Sigma\vdash \phi$, a contradiction.

\end{proof}
\end{Theorem}

\subsection{Completeness Theorem for $\mathfrak{M}$-Valued Logic}

\begin{Definition}
Let $\Sigma$ be a set of formulas. $\phi \equiv_\Sigma \theta$ if and only if $\Sigma \vdash \phi  \Leftrightarrow \theta$.
\end{Definition}

\begin{Definition}
Let $F$ be the set of all formulas. $\mathfrak{F}=\langle F,\vee,\wedge, \neg, t_p, \exists v_k, v_k\thickapprox v_l,\gamma_p \rangle_{k,l\in \omega, p\in M}$
\end{Definition}
 $\mathfrak{F}$ is an algebra.
 
 \begin{Theorem}
$\equiv_\Sigma$ is a congruence relation on  $\mathfrak{F}$.
\begin{proof}
Assume that $\phi_1 \equiv_\Sigma \theta_1$ and $\phi_2 \equiv_\Sigma \theta_2$. To show that $\exists v_k\phi_1 \equiv_ \Sigma\exists v_k \theta_1$, use Theorem \ref{Theorems} (c). To show that  $\phi_1\wedge\phi_1 \equiv_\Sigma \phi_2\wedge\phi_2$, use the tautology $(\phi_1\Leftrightarrow\theta_1)\rightarrow((\phi_2\Leftrightarrow\theta_2) \rightarrow(\phi_1\wedge\phi_2)\Leftrightarrow(\theta_1\wedge\theta_2))$. Similarily, the rest of the proof is done using appropriate tautologies.
\end{proof}
\end{Theorem}

Denote the equivalence classes in $\mathfrak{Fm}/\equiv_\Sigma$  by $[\phi]$. It is easily seen that $s\binom{k}{l}[\phi]=[\boldsymbol{S}\binom{k}{l}\phi]$
 \begin{Theorem}
 $\mathfrak{F}/ \equiv_\Sigma$ is a $\mathfrak{M}-\textrm{CA}_{\omega}$.
 \begin{proof}
 To verify all the axioms of $\mathfrak{M}-\textrm{CA}_{\omega}$s, except for Axiom 30, use tautologies  22 to 24 and 27 to 46 of Theorem \ref{Tautologies} and validities 2 to 13 (which are axioms of our logic) of Theorem \ref{validities}. To verify Axiom 30, use Theorem \ref{Theorems} (f), Theorem \ref{Tautologies} 47 and Modus Ponens.
 \end{proof}
\end{Theorem}

In $\mathfrak{Fm}/\equiv_\Sigma$, $dim([\phi])=\bigcup_{p\in M}\{k:c_k\delta_p[\phi]\ne\delta_p[\phi]\}=\bigcup_{p\in M}\{k:[\exists v_k\gamma_p\phi]\ne[\gamma_p\phi]\}$.
 
 \begin{Theorem}
 $\mathfrak{Fm}/\equiv_\Sigma$ is a locally finite $\mathfrak{M}-\textrm{CA}_{\alpha}$.
 \begin{proof}
 Let $k\notin Fv(\phi)$. By Theorem \ref{Theorems} (e), $\Sigma\vdash\phi\Leftrightarrow\exists v_k\phi$ and hence, $[\phi]\Leftrightarrow[\exists v_k\phi]=c_k[\phi]$. Since $Fv(\gamma_p\phi)=Fv(\phi)$, $k\notin Fv(\gamma_p\phi)$ and $\Sigma\vdash\gamma_p\phi\Leftrightarrow\exists v_k\gamma_p\phi$. Therefore, $[\gamma_p\phi]=[\exists v_k\gamma_p\phi]$ for all $p\in M$ and $k\notin dim([\phi])$. Therefore, $dim([\phi)\subseteq Fv(\phi)$. Since $Fv(\phi)$ is finite, so is $dim([\phi])$.
 
 \end{proof}
\end{Theorem}

\begin{Theorem}
\label{more than 1}
If $\Sigma$ is consistent, then $|\mathfrak{Fm}/\equiv_\Sigma|$ has more than one element.
\begin{proof}
Suppose that $\Sigma$ is consistent, then $|\mathfrak{Fm}/\equiv_\Sigma|$ has only one element. Then $[\theta]=[\phi]$ for all $\theta $ and $\phi$. Therefore, $[t_1]=[t_0]$ and hence, $t_1\equiv_\Sigma t_0$ and  $\Sigma \vdash t_1 \Leftrightarrow t_0$. Use the tautology $(t_1\Leftrightarrow t_0)\rightarrow t_0$ and Modus Ponens to obtain $\Sigma \vdash t_0$ and hence $\Sigma$ is inconsistant.
\end{proof}
\end{Theorem}

\begin{Theorem}
If $\phi\in \Sigma$, then $\Sigma\vdash\phi\Leftrightarrow t_1$.
\begin{proof}
Here is a proof: $\phi$, $\Gamma\phi$, $\Gamma\phi\rightarrow(\phi \rightarrow(\phi\Leftrightarrow t_1))$ (tautology) $\phi \rightarrow(\phi\Leftrightarrow t_1)$, $\phi\Leftrightarrow t_1$
\end{proof}
\end{Theorem}

\begin{Theorem}(Completeness theorem 1)\\
(a) Every consistent set of formulas has a model.\\
(b) Every $Q$-consistent set of formulas has a $Q$-model where $Q\subseteq M$.
\begin{proof}
(a) Let $\Sigma$ be a consistent set of formulas. $\mathfrak{Fm}/\equiv_\Sigma$ is a locally finite $\mathfrak{M}-\textrm{CA}_{\omega}$. By Theorem \ref{rep 2}, $\mathfrak{Fm}/\equiv_\Sigma$ is isomorphic to a subalgebra of a product of regular $\mathfrak{M}-\textrm{CSA}_{\omega}$'s, $\mathfrak{C}_i$ for $i \in J$.
Let $f:\mathfrak{Fm}/\equiv_\Sigma\rightarrow\prod_{i\in J}\mathfrak{C}_i$ be an injection and let $P_j:\prod_{i\in j}\mathfrak{C}_i\rightarrow\mathfrak{C}_j$ be the projection function. Since $\Sigma$ is consistent, by Theorem \ref{more than 1}, for some $j\in J$, $|\mathfrak{C_j}|$ has more than one element. Let $h=P_j\circ f$. Then $h:\mathfrak{Fm}/\equiv_\Sigma\rightarrow\mathfrak{C_j}$ is a homomorphism.\\
It is easy to show that $s\binom{k}{l}[\phi]=[S\binom{k}{l}\phi]$ and $h(s\binom{k}{l}[\phi])=\boldsymbol{S}\binom{k}{l}h[\phi]$.

Let $\mathfrak{A}=\langle |\mathfrak{C}_j|,P_i\rangle_{i\in I}$ where $P_i$ is the $\mathfrak{M}$-valued relation on $|\mathfrak{C}_j|$ defined by

\begin{center}
$(P_i)^q(x_0,\cdots ,x_{n_i})$ if and only if $t\binom{0}{x_0}\cdots \binom{n_i-1}{x_{n_i-1}}\in h([R_iv_0\cdots v_{n_i-1}]^q\}$
\end{center}

Note that this definition is independent of $t$ since $\mathfrak{C_j}$ is regular. The $\mathfrak{M}-\textrm{CSA}_{\omega}$ $\mathfrak{S}$ of Definition 3.5 is a subalgebra of $\mathfrak{C_j}$.

$(R_iv_0,\cdots,v_{n_i-1})^\mathfrak{A}=h([R_i{v_0},\cdots,v_{n_{i-1}}])$ since

\begin{eqnarray*}
t\in ((R_iv_0,\cdots,v_{n_{i-1}})^\mathfrak{A})^q & \textrm{if and only if}& (P_i)^q(t_0,\cdots,t_{n_{i-1}})(\textrm{by definition of } \phi^\mathfrak{A})\\
&\textrm{if and only if}& t\binom{0}{t_0}\cdots\binom{n-1}{t_{n_{i-1}}}\in h([R_i{v_0},\cdots,v_{n_{i-1}}])^q(\textrm{by definition of } P_i)\\
&\textrm{if and only if}& t\in h([R_{v_0},\cdots,v_{n_{i-1}}])^q\\
\end{eqnarray*}
Remember that in the $\mathfrak{M}-\textrm{CSA}_{\alpha}$, $\mathfrak{S}$, $(S^k_l\phi)^\mathfrak{A}=\boldsymbol{S}^k_l\phi^\mathfrak{A}$.
We show by induction on formulas that $h([\phi])=\phi^\mathfrak{A}$.\\
$\phi=R_iv_{j_0}\cdots,v_{n_{j_{i-1}}}$:\\
Let $k_m\ne k_l$ for $m\ne l$ and  $k_m\notin \{0,\cdots n-1,j_0\cdots j_{n_{i-1}}\}$ for $m<n_i$\\

\begin{eqnarray*}
h([R_iv_{j_0}\cdots,v_{n_{j_{i-1}}}])&=&h([S^{k_0}_{j_0}\cdots S^{k_{n_i-1}}_{j_{n_i-1}}S^0_{k_0}\cdots S^{n_{i-1}}_{k_{n_i-1}}R_iv_0\cdots v_{n_{i-1}}]) \\
&=&h(s^{k_0}_{j_0}\cdots s^{k_{n_i-1}}_{j_{n_i-1}}s^0_{k_0}\cdots s^{n_{i-1}}_{k_{n_i-1}}[R_iv_0\cdots v_{n_{i-1}}])\\
&=&\boldsymbol{S}^{k_0}_{j_0}\cdots \boldsymbol{S}^{k_{n_i-1}}_{j_{n_i-1}}\boldsymbol{S}^0_{k_0}\cdots \boldsymbol{S}^{n_{i-1}}_{k_{n_i-1}}h([R_iv_0\cdots v_{n_{i-1}}])\\
&=&\boldsymbol{S}^{k_0}_{j_0}\cdots \boldsymbol{S}^{k_{n_i-1}}_{j_{n_i-1}}\boldsymbol{S}^0_{k_0}\cdots \boldsymbol{S}^{n_{i-1}}_{k_{n_i-1}}(R_iv_0,\cdots,v_{n_i-1})^\mathfrak{A}\\
&=&(S^{k_0}_{j_0}\cdots S^{k_{n_i-1}}_{j_{n_i-1}}S^0_{k_0}\cdots S^{n_{i-1}}_{k_{n_i-1}}R_iv_0,\cdots,v_{n_i-1})^\mathfrak{A}\\
&=&(R_iv_{j_0}\cdots,v_{n_{j_{i-1}}})^\mathfrak{A}\\
\end{eqnarray*}
\end{proof}
\end{Theorem}
The first and last steps follow from 14 of Theorem \ref{validities} which is an axiom.\\

$h([t_p])=\boldsymbol{U}_p=t_p^\mathfrak{A}$ in the $\mathfrak{M}-\textrm{CSA}_{\alpha}$ $\mathfrak{C_j}$ or $\mathfrak{S}$ since $h:\mathfrak{Fm}/\equiv_\Sigma\rightarrow\mathfrak{C_j}$  is a homomorphism.\\
$h([v_n \thickapprox v_m])=\boldsymbol{D}_{ij}$ in the $\mathfrak{M}-\textrm{CSA}_{\alpha}$ $\mathfrak{C_j}$ or $\mathfrak{S}$ since $h:\mathfrak{Fm}/\equiv_\Sigma\rightarrow\mathfrak{C_j}$  is a homomorphism.\\
Now assume that $h([\phi])=\phi^\mathfrak{A}$ and $h([\theta])=\theta^\mathfrak{A}$.\\
$h([\neg\phi])=h(-[\phi])=\circleddash h([\phi])= \circleddash\phi^\mathfrak{A}=(\neg\phi)^\mathfrak{A}$\\
$h([\phi\wedge\theta])=h([\phi]\cdot[\theta])=h([\phi])\sqcap h([\theta])= \phi^ \mathfrak{A} \sqcap\theta^ \mathfrak{A} =(\phi\wedge\theta)^ \mathfrak{A} $\\
$h([\phi\vee\theta])=h([\phi]+[\theta])=h([\phi])\sqcup h([\theta])= \phi^ \mathfrak{A} \sqcup\theta^ \mathfrak{A} =(\phi\vee\theta)^ \mathfrak{A} $\\
$h([\exists v_k\phi])=h([c_k[\phi])=\boldsymbol{C}_kh([\phi])=\boldsymbol{C}_k\phi^ \mathfrak{A}=(\exists v_k\phi)^ \mathfrak{A}$\\
$h([\gamma_p\phi])=h(\delta_p[\phi])=\boldsymbol{\delta}_ph([\phi])=\boldsymbol{\delta}_p\phi^ \mathfrak{A}=(\gamma_p\phi)^ \mathfrak{A}$\\
Let $\phi \in \Sigma$. Then $[\phi]=[t_1]$ in the $\mathfrak{Fm}/\equiv_\Sigma$ Then $\phi^\mathfrak{A}=h([\phi])=h([t_1])=\boldsymbol{U}_1$ in $\mathfrak{S}$. Therefore, $\phi$ is true in $\mathfrak{A}$ and $\mathfrak{A}$ is a model of $\Sigma$.\\
(b) Let $\Sigma$ be $Q$-consistent. Therefore, by Theorem \ref{Q-consistent}, $\Sigma^Q$ is consistent and by part (a), has a model $\mathfrak{B}$. Therefore, by Theorem \ref{Q model}, $\mathfrak{B}$ is a $Q$-model of $\Sigma$.

\begin{Theorem}
\label{sentence case}
If $\Sigma$ is set of formulas and $\phi$ is a sentence, then $\Sigma\models\phi$ implies $\Sigma\vdash\phi$.
\begin{proof}
Assume that $\Sigma\models\phi$. First we show that $\Sigma\vdash\Gamma\phi$. Suppose that it is not the case that $\Sigma\vdash\Gamma\phi$. Since $\Gamma\Gamma\phi$ is a tautology, we have $\Sigma\vdash \Gamma\Gamma\phi$. By Theorem \ref{consistency}, $\Sigma\cup\{\neg \Gamma\phi\}$ is consistent and hence has a model $\mathfrak{A}$. Therefore, $\mathfrak{A}$ is a model of $\Sigma$ and $\neg \Gamma\phi$ is true in $\mathfrak{A}$. Since $\Sigma\models\phi$ we also have $\Sigma\models\Gamma\phi$. Therefore, $\Gamma\phi$ is true in $\mathfrak{A}$ which is a contradiction and we have $\Sigma\vdash\Gamma\phi$.\\

Now assume that it is not the case that $\Sigma\vdash\phi$. By Theorem \ref{consistency}, $\Sigma\cup\{\neg\phi\}$ and has a model $\mathfrak{B}$. $\mathfrak{B}$ is a model of $\Sigma$ and $\neg\phi$ is true in $\mathfrak{B}$. Since $\Sigma\models\phi$, $\phi$ is true in $\mathfrak{B}$ which is a contradiction. Therefore, $\Sigma\vdash\phi$.
\end{proof}
\end{Theorem}

\begin{Theorem} (Completeness Theorem 2)\\
Let $Q\subseteq M$ Then\\
(a) If $\Sigma$ is set of formulas then $\Sigma\models\phi$ implies $\Sigma\vdash\phi$.\\
(b) If $\Sigma$ is set of formulas then $\Sigma\models_Q\phi$ implies $\Sigma\vdash_Q\phi$.
\begin{proof}
(a) Assume that $\Sigma\models\phi$ and let $\theta=\forall v_{k_0}\cdots\forall v_{k_{n-1}}\phi$ where $Fv(\phi)=\{k_0,\cdots,{k_{n-1}}\}$. Then $\theta$ is a sentence. Let $\mathfrak{A}$ be a model of $\Sigma$. Then $\phi$ is true in $\mathfrak{A}$ and $\phi^\mathfrak{A}=\boldsymbol{U}_1$ in the $\mathfrak{M}-\textrm{CSA}_{\alpha}$ $\mathfrak{S}$. It is easy to verify that $\boldsymbol{Q}_{m}\boldsymbol{U}_1=\boldsymbol{U}_1$. $\theta^\mathfrak{A}=(\forall v_{k_0}\cdots\forall v_{k_{n-1}}\phi)^\mathfrak{A}=\boldsymbol{Q}_{k_0}\cdots\boldsymbol{Q}_{k_{n-1}}\phi^\mathfrak{A}=\boldsymbol{Q}_{k_0}\cdots\boldsymbol{Q}_{k_{n-1}}\boldsymbol{U}_1=\boldsymbol{U}_1$. Therefore, $\Sigma\models\theta$ and by part Theorem \ref{sentence case}, $\Sigma\vdash\theta$.\\
Let $\psi_1,\cdots,\psi_m=\forall v_{k_0}\cdots\forall v_{k_{n-1}}\phi$ be a proof of $\Sigma\vdash\forall v_{k_0}\cdots\forall v_{k_{n-1}}\phi$. Then $\psi_1,\cdots,\psi_m,\Gamma\forall v_{k_0}\cdots\forall v_{k_{n-1}}\phi,$\\
$\Gamma\forall v_{k_0}\cdots\forall v_{k_{n-1}}\phi\rightarrow(\forall v_{k_0}\cdots\forall v_{k_{n-1}}\phi\rightarrow\forall v_{k_1}\cdots\forall v_{k_{n-1}} \phi)\textrm{ (axiom)},\forall v_{k_0}\cdots\forall v_{k_{n-1}}\phi\rightarrow\forall v_{k_1}\cdots\forall v_{k_{n-1}} \phi),\forall v_{k_1}\cdots\forall v_{k_{n-1}} \phi)$ is a proof of $\Sigma\vdash\forall v_{k_1}\cdots\forall v_{k_{n-1}}\phi$\\
Continuing this way yields $\Sigma\vdash\phi$\\
(b) Let $\Sigma\models_Q\phi$. By Theorem \ref{Q models}, $\Sigma^Q\models\phi^Q$. By Part (a), $\Sigma^Q\vdash\phi^Q$ and, by definition, $\Sigma\vdash_Q\phi$.
\end{proof}
\end{Theorem}

\begin{center}
References
\end{center}

[1] L. Bolc and P. Borowick, Many valued-logics, vol. 1, Springer-Verlag, Berlin, 1992.

[2] N. \textsc{Feldman}, \textit{The cylindric algebras of three-valued logic}, \textit{\textbf{Journal of symbolic logic}}, vol. 63 (1998), pp. 1201-1217.

[3] L.\textsc{Henkin}, J.D.  \textsc{Monk}, and A. \textsc{Tarski},  \textit{\textbf{Cylindric algebras, part I}}, North-Holland Publishing Company, Amsterdam,1971.

[4] L.\textsc{Henkin}, J.D.  \textsc{Monk}, and A. \textsc{Tarski},  \textit{\textbf{Cylindric algebras, part II}}, North-Holland Publishing Company, Amsterdam,1985.

\end{document}